\documentclass{amsart}
\usepackage{amscd}
\usepackage[dvips]{graphicx}
\usepackage{amscd}
\usepackage{amsmath}
\usepackage{amsfonts}
\usepackage{amssymb}
\usepackage[latin1]{inputenc}
\newtheorem{theo+}    {Theorem}      [section]
\newtheorem{prop+}  [theo+]  {Proposition}
\newtheorem{coro+}  [theo+]  {Corollary}
\newtheorem{lemm+}  [theo+]  {Lemma}
\newtheorem{deep+}  [theo+]  {Deep Result}
\newtheorem{fact+}  [theo+]  {Fact}
\theoremstyle{definition}
\newtheorem{exam+}  [theo+]  {Example}
\newtheorem{rema+}  [theo+]  {Remark}
\newtheorem{defi+}  [theo+]  {Definition}
\newtheorem{xca+}[theo+]{Exercise}
\newenvironment{theorem}{\begin{theo+}}{\end{theo+}}
\newenvironment{proposition}{\begin{prop+}}{\end{prop+}}
\newenvironment{corollary}{\begin{coro+}}{\end{coro+}}
\newenvironment{lemma}{\begin{lemm+}}{\end{lemm+}}

\newenvironment{example}{\begin{exam+}}{\end{exam+}}
\newenvironment{remark}{\begin{rema+}}{\end{rema+}}
\newenvironment{definition}{\begin{defi+}}{\end{defi+}}

\numberwithin{equation}{section}

\def\beqn{\begin{equation}}
\def\eeqn{\end{equation}}

\def\epf{\qed \enddemo}

\def\a{\alpha}

\def\Claminv2{|C(\Lambda)|^{-2}}

\def\lam{\lambda}

\def\Ome{\Omega}

\def\Aa2D{A^{\a,2}(D)}
\def\bAa2D{\overline{A^{\a,2}(D)}}
\def\Ab2D{A^{\beta,2}(D)}
\def\bAb2D{\overline{A^{\beta,2}(D)}}

\def\Norm#1_#2{\Vert#1\Vert_{#2}}

\def\2pd#1#2{\frac{\partial^2 #1}{\partial #2^2}}
\def\p11d#1#2#3{\frac{\partial^2 #1}{  \partial #2\partial #3  }}

\def\Claminv2{|C(\Lambda)|^{-2}}

\def\lam{\lambda}

\def\Aa2D{A^{\a,2}(D)}
\def\bAa2D{\overline{A^{\a,2}(D)}}
\def\Ab2D{A^{\beta,2}(D)}
\def\bAb2D{\overline{A^{\beta,2}(D)}}

\def\ub1#1{\underline{\mathbf 1^{#1}}}

\def\nat0{\mathbb Z_{\ge 0}}

\def\bpf{\begin{proof}}
\def\epf{\end{proof}}
\def\beq{\begin{equation}}
\def\eeq{\end{equation}}

\def\draft{\centerline{(Draft {\the \day}/{\the\month} \the \year.)}}

\begin{document}
\title[Segal-Bargmann transform]{
Segal-Bargmann  transform
 on Hermitian symmetric spaces and Orthogonal Polynomials}
\author{Mark Davidson, Gestur \'{O}lafsson, and Genkai Zhang}
\address{Department of Mathematics, Louisiana State University, Baton Rouge, LA\ 70803,
USA}
\address{Department of Mathematics, Louisiana State University, Baton Rouge, LA\ 70803,
USA} \address{Department of Mathematics, Chalmers University of
Technology and G\"oteborg University, S-412 96 Got\"oborg, Sweden}
\email{davidson@math.lsu.edu} \email{olafsson@math.lsu.edu}
\email{genkai@math.chalmers.se}
\thanks{Research by G. \'{O}lafsson supported by NSF grant DMS 0070607 and the MSRI}
\thanks{Research by G. Zhang supported by Swedish Research Council (VR)}
\begin{abstract}Let $\mathcal{D}=G/K$ be a complex bounded symmetric domain of tube type in a
complex Jordan algebra $V$ and let
$\mathcal{D}_{\mathbb{R}}=H/L\subset \mathcal{D}$ be its real form
in a formally real Euclidean Jordan algebra $J\subset V$. We
consider representations of $H$ that are gotten by the generalized
Segal-Bargmann transform from a unitary $G$-space of holomorphic
functions on $\mathcal{D}$ to an $L^2$-space on
$\mathcal{D_{\mathbf{R}}}$.  We prove that in the unbounded
realization the inverse of the unitary part of the restriction map
is actually the Laplace transform. We find the extension to
$\mathcal{D}$ of the spherical functions on
$\mathcal{D}_{\mathbb{R}}$ and find the expansion in terms of the
$L$-spherical polynomials on $\mathcal{D}$, which are Jack
symmetric polynomials.  We prove that the coefficients are
orthogonal polynomials in an $L^2$-space, the measure being the
Harish-Chandra Plancherel measure multiplied by the symbol of the
Berezin transform. We prove the difference
equation and recurrence relation for  those polynomials by considering
the action of the Lie algebra  and the Cayley transform
on the polynomials on $\mathcal D$.
\end{abstract}
\maketitle

\section{Introduction}

The study of various generalizations of the classical Weyl
transform has attracted much interest and has been pursued for
some time. As is well-known the Weyl transform maps unitarily from
the $L^{2}$-space on $\mathbb{C}^{n}=\mathbb{R}^{2n}$ onto the
space of Hilbert-Schmidt operators on Fock-space,
$\mathcal{F}(\mathbb{C}^{n})$, on $\mathbb{C}^{n}$. The Weyl
transform intertwines the natural actions of the Heisenberg group,
$H_{n}=\mathbb{R}\times\mathbb{C}^{n}$ and realizes the
decomposition of $L^2(H_n)$.  However, there is  another model for
the representation, namely, the so-called Schr\"{o}dinger model,
and the Segal-Bargmann transform is then a unitary map from the
Schr\"{o}dinger model to the Fock model intertwining the action of
the Heisenberg group. The Segal-Bargmann and Weyl transforms have
been studied for a long time in its connection with geometric
quantization. But somewhat surprisingly their unitary properties
have not been satisfactorily  clarified. The unifying idea is the
\textit{restriction principle}, i.e. polarization of a suitable
restriction map(\cite{OOe96},\cite{oz-weyl}). This idea, in turn,
also unifies the so-called Wick quantization and Berezin transform
in one picture, which we briefly recall.

Consider the tensor product
$\mathcal{F}(\mathbb{C}^{n})\otimes\overline
{\mathcal{F}(\mathbb{C}^{n})}$ realized as the space of
Hilbert-Schmidt operators with integral kernels $F(z,w)$
holomorphic in $z$ and
anti-holomorphic in $w$. Imbed $\mathbb{C}^{n}$ into $\mathbb{C}^{n}%
\times\overline{\mathbb{C}}^{n}$, bar denoting opposite complex structure, by
$z\mapsto(z,z)$. Consider the restriction mapping from $\mathcal{F}%
(\mathbb{C}^{n})\otimes\overline{\mathcal{F}(\mathbb{C}^{n})}$ to real
analytic functions on $\mathbb{C}^{n}$ taking functions $F(z,w)$ to its
restriction $F(z,z)$. By taking a multiplier of the Gaussian (which is the
restriction of the reproducing kernel of $\mathcal{F}(\mathbb{C}^{n})$) into
account one gets a bounded injective map $R:\mathcal{F}(\mathbb{C}%
^{n})\otimes\overline{\mathcal{F}(\mathbb{C}^{n})}\rightarrow L^{2}%
(\mathbb{C}^{n})$ with dense image. It turns out that the adjoint
$R^{\ast}$ is the Wick quantization map. Consider the polar
decomposition $R^{\ast}=U \sqrt{R R^{\ast}}$. The map $U$ is the
Weyl quantization and  $RR^{\ast}$ is the Berezin transform. Thus
the information about the Wick quantization, the Weyl
quantization, and the Berezin transform are encoded in the
restriction map.

It is easy to see by a similar calculation, as in \cite{oz-weyl},
that by restriction of the Fock space on $\mathbb{C}^{n}$ to its
real form and taking the unitary part of the restriction we get
the Bargmann-Segal transform. The Weyl transform from
$L^{2}(\mathbb{C}^{n})=L^{2}(\mathbb{R}^{n})\otimes
L^{2}(\mathbb{R}^{n})$ onto $\mathcal{F}(\mathbb{C}^{n})\otimes
\overline{\mathcal{F}(\mathbb{C}^{n})}$ is the tensor product of
two Bargmann-Segal transforms, obtained by considering the
restriction to two different real forms. The point is then that
both transforms may be considered as the unitary part of the
adjoint of the  restriction of certain holomorphic representations
on various real forms of the underlying complex manifolds.

Instead of a flat complex space $\mathbb{C}^{n}$ we may take a
bounded symmetric domain $\mathcal{D}=G/K$ and consider the tensor
product of a weighted Bergman space (so called holomorphic
discrete series of $G$) with its conjugate. One may then perform
the polar decomposition $R^{\ast}=U \sqrt{RR^{\ast}}$ of the
restriction map $R$ and get a unitary intertwining operator $U$
from the tensor product onto $L^{2}(\mathcal{D})$. In terms of
representation language this has been studied by Repka
\cite{Repka}, and an analytic and explicit
approach was started
in \cite{gkz-invtoe}. In particular, it is realized  that
the analytic issues are far more subtle
when one considers the tensor product of the analytic continuation
of weighted Bergman spaces (\cite{gkz-tenmin}). For that purpose
we need to understand the positive part in the polar
decomposition, namely the square root of the Berezin transform
$RR^{\ast}$. The Berezin transform in the case
of weighted Bergman spaces is a positive
bounded operator on $L^{2}(\mathcal{D})$ and its symbol has been
calculated by Unterberger-Upmeier \cite{UU}.

In \cite{OOe96}, it is shown how to generalize those ideas to
restriction maps from a reproducing Hilbert space of holomorphic
function on a complex manifold $M_{\mathbb{C}}$ to a totally real
submanifold $M$. A particular situation is when
$M_{\mathbb{C}}=G/K$ is a bounded symmetric domain, the Hilbert
space $\mathcal{F}(\mathbb{C}^{n})$ is replaced by a weighted
Bergman space, and $M=H/H\cap K$ is a totally real homogeneous
submanifold. The restriction principle gives a natural way to
define a Segal-Bargmann and Berezin transform. The symbol of the
Berezin transform is then calculated in \cite{gkz-bere-rbsd}. The
same result for real tube domains and for classical domains are
obtained by van Dijk and Pevner \cite{vanDijk-Pevsner} and,
respectively, Neretin \cite{Neretin}.

In this paper we will use these ideas to construct and study a
natural class of functions and orthogonal polynomials on real
symmetric domains.  We consider the simplest case when the real
symmetric domain is, in the Siegel domain realization, the
symmetric cone (Type A in terms of the classification of the root
system, see \cite{gkz-bere-rbsd}). The starting point is a unitary
highest weight representation $(\pi_{\nu},\mathcal{H}_{\nu})$,
with one dimensional minimal $K$-type, of the group $G$. We
consider the restriction map to a real form $H/H\cap K$
corresponding to the symmetric cone, both for the bounded
realization and the unbounded realization of $G/K$. The main
difference is, that in the bounded realization we include the
multiplier corresponding to the minimal $K$-type in our
restriction map
\[
R:\mathcal{H}_{\nu}\rightarrow L^{2}(H/H\cap K),\qquad f\mapsto
D_{\nu}f\vert_{H/H\cap K},
\]
whereas in the unbounded realization we include this multiplier in
the measure, so that the map $R:\mathcal{H}_{\nu}\rightarrow
L^{2}(\Omega ,d\mu_{\nu})$ is now just the restriction $f\mapsto
f|_{\Omega}$. Then in the unbounded realization the Segal-Bargmann
transform turns out to be the Laplace transform (see Theorem
\ref{th-Laplace}). Considering the unitary image of an orthogonal
set of $K$-finite $H\cap K$-invariant vectors in the holomorphic
realization we construct a family of spherical orthogonal
functions and polynomials on $H/H\cap K$. Our main result is the
orthogonality relations, recurrence formulas and difference
formulas for the polynomials. In the unbounded realization we get
functions of the form $\omega\mapsto
e^{-\text{Tr}(\omega)}L_{\nu}(2\omega)$ where $L_{\nu}$ is a
polynomial which agrees with the Laguerre polynomials in the case
where $G=\mathrm{SL}(2,\mathbb{R})$ and $H/H\cap
K\simeq\mathbb{R}^{+}$ (\cite{doz1}). These polynomials have also
been considered in \cite{FK-book}, Chapter XV, but here we relate
them to representations of $G$. In that way we derive differential
equations satisfied by these polynomials (see Theorem
(\ref{The-Euler})). Furthermore, applying the spherical Fourier
transform we get an orthonormal basis of
$L^{2}(\mathfrak{a}^{\ast},\left| c(\lambda)\right|
^{-2}d\lambda)^{W}$ and a family of orthogonal polynomials on $\mathfrak{a}%
^{\ast}$ (see Corollary(\ref{cor-orthog}) and Proposition
(\ref{prop-poly})). We derive recurrence formulas and difference
formulas for these polynomials (see Theorems (\ref{thm-recrel})
and (\ref{th-diff})).

We mention that a large class of symmetric and non-symmetric
polynomials satisfying various difference equations has been
recently introduced by using the representation theory of affine
Hecke algebras. However, as it is noted by Cherednik
(\cite{Cherednik-imrn-95}, p. 484) the meaning of the difference
equations needs to be clarified. Our results provide the
clarification and meaning of the orthogonality and difference
relations  for these orthogonal polynomials. It also reveals that
the difference relations and the recurrence relations are somewhat
dual to each other. We may well expect that this will be true for
other types of Macdonald polynomials.

This joint project started in December 2000 as G. Zhang was
visiting the Louisiana State University. In our discussions we
realized that we were working on similar problems using similar
ideas, except that G. Zhang was working in the bounded realization
for all real tube domains, while the first two authors were
working on the restriction to the symmetric cone (Type A case) of
unbounded realizations of Siegel domains. The analytic nature
involved for the Type A case is somewhat richer, in particular in
its connection with the Laplace transform. The treatment of types
C and D requires however some different and more combinatorial
methods, see \cite{gkz-jfa}.

\subsection*{Acknowledgement}

We would like to thank Hjalmar Rosengren for bringing the paper
\cite{Basu-Wolf} to our attention. The third author would also
like to thank NSF for partially supporting his visit to LSU\ in
December 2000.

\section{Bounded symmetric domains, symmetric cones and Cayley transform}

In this section we recall some known necessary background about
bounded symmetric domains. We follow the presentation in
\cite{Loos-bsd}; see also \cite{FK-book}.

\subsection{Bounded symmetric domains}

Let $V$ be a $d$-dimensional \textit{complex} Hermitian simple
Jordan algebra. Let $\mathcal{D}\subset V$ be an irreducible
bounded symmetric domain isomorphic to a tube type domain in $V$.
Let $\mathrm{Aut}(\mathcal{D})$ be the group of all biholomorphic
automorphisms of $\mathcal{D}$, let
$G=\mathrm{Aut}(\mathcal{D})_{0}$ be the connected component of
the identity in $\mathrm{Aut}(\mathcal{D})$, and let $K$ be the
isotropy subgroup of $G$ at the point $0$. Then $G$ is semisimple
and, as a Hermitian symmetric space, $\mathcal{D}=G/K$.
Furthermore, $K$ is a maximal compact subgroup of $G$. We denote
by \ $\tilde{G}$ the connected simply connected covering group of
$G$
and $\tilde{K}$ the pre-image of $K$ in $\tilde{G}$. Then $\tilde{K}%
\simeq\mathbb{R}\times K_{1}$ where $K_{1}$ is a simply connected
compact semisimple Lie group. We start by reviewing the basic
structure theory of $G$ in terms of the structure of the Jordan
algebra $V$. The Lie algebra
$\mathfrak g$ of $G$ is identified with the Lie algebra $\mathrm{aut}%
(\mathcal{D})$ of all completely integrable holomorphic vector
fields on $\mathcal{D}$ equipped with the Lie product
\[
\lbrack X,Y](z):=X^{\prime}(z)Y(z)-Y^{\prime}(z)X(z),\,\,\,X,Y\in
\mathrm{aut}(\mathcal{D}),\,z\in\mathcal{D}.
\]

Define $\theta:\mathfrak{g}\rightarrow\mathfrak{g}$ by
$\theta(X)(z):=-X(-z)$. Then $\theta$ is a Cartan involution on
$\mathfrak{g}$ and
\[
\mathfrak{g}^{\theta}=\left\{
X\in\mathfrak{g}\mid\theta(X)=X\right\} =\mathfrak{k}\,,
\]
is the Lie algebra of $K$. Let
\[
\mathfrak{p}:=\{X\in\mathfrak{g}\mid\theta(X)=-X\}\,.
\]
Then $\mathfrak{g}=\mathfrak{k}\oplus\mathfrak{p}$ is the Cartan
decomposition of $\mathfrak{g}$ corresponding to $\theta$. Every
element in $\mathfrak{g}$ is given by a polynomial of degree at
most $2$, \cite{UP87}. The Lie algebra $\mathfrak{k}_{\mathbb C}$
corresponds to the elements of degree one and is isomorphic to a
subalgebra of $\mathrm{End}(V)$ by
$T\mapsto(Tz)\frac{\partial\,}{\partial z}$. The identity map thus
corresponds to the \textit{Euler operator}
\begin{equation}
Z_{0}:=z\frac{\partial\,}{\partial z}\,. \label{eq-Z0}%
\end{equation}
This element is central in $\mathfrak{k}_{\mathbb{C}}$ and
$\mathrm{ad}(Z_{0})$ has eigenvalues $\pm1$ on
$\mathfrak{p}_{\mathbb{C}}$. The $+1$-eigenspace
$\mathfrak{p}^{+}$
corresponds to the constant polynomials and the $-1$-eigenspace $\mathfrak{p}%
^{-}$ corresponds to the polynomials of degree two. There exists a quadratic
form $Q:V\rightarrow\mathrm{End}(\bar{V},V)$ (where $\overline{V}$ is the
space $V$ but with the opposite complex structure), such that
\begin{equation}
\mathfrak{p}=\{\xi_{v}\mid v\in
V\},\quad\mathrm{where}\quad\xi_{v}(z):=(v-Q(z)\bar
{v})\frac{\partial\,}{\partial z}\,. \label{eq-p}%
\end{equation}

In the following we will identify elements in $\mathfrak{g}$ with
the corresponding polynomials. Let $\{z,\bar{v},w\}$ be the
polarization of $Q(z)\bar{v}$, i.e.,
\[
\{z,\bar{v},w\}=(Q(z+w)-Q(z)-Q(w))\bar{v}.
\]
Then
\begin{equation}
\xi_{v}(z)=v-Q(z)\bar{v}=v-\frac{1}{2}\{z,\bar{v},z\}. \label{eq-xiv}%
\end{equation}

The space $V$ with this triple product on
$V\times\overline{V}\times V$, is a JB*-triple; (\cite{UP87}).
Define $D:V\times\bar{V}\rightarrow \mathrm{End}(V)$ by
$D(z,\bar{v})w=\{z,\bar{v},w\}$. The Lie bracket of two elements
$\xi_{v},\xi_{w}\in\mathfrak{p}$ is then given by
\begin{equation}
\lbrack\xi_{v},\xi_{w}]=(D(v,\bar{w})-D(w,\bar{v}))\in\mathfrak{k}\,.
\label{LieInP}%
\end{equation}
The group $K$ acts on $V$ by linear transformations. Furthermore,
\begin{equation}
z,w\mapsto(z\mid w):=\frac{1}{p}\mathrm{Tr}D(z,\bar{w}) \label{eq-innerpr}%
\end{equation}
is a $K$-invariant inner product on $V$. Here $\mathrm{Tr}$ is the
trace functional on $\mathrm{End}(V)$ and $p=p(\mathcal{D})$ is
the genus of $\mathcal{D}$ (see (\ref{eq-gen}) below).  Denote by
$||\cdot||$ the corresponding $K$-invariant operator norm on
$\mathrm{End}(V)$. Define the \textit{spectral norm} $||\cdot||_{\mathrm{sp}}$  on $V$ by
\[
\Vert z\Vert_{\mathrm{sp}}:=\Vert\frac{1}{2}D(z,\bar{z})\Vert^{1/2}.
\]
Then the domain $\mathcal{D}$ is realized as the open unit ball of $V$ with
respect to the spectral norm, i.e. $\mathcal{D}=\{z\in V\mid\Vert
z\Vert_{\mathrm{sp}}<1\}$.

An element $v\in V$ is called a \textit{tripotent} if
$\{v,\bar{v},v\}=2v$. Let us choose and fix a frame
$\{e_{j}\}_{j=1}^{r}$ of tripotents in $V$. The
number $r$ is called the \textit{rank} of $\mathcal{D}$. Define $H_{j}%
=D(e_{j},\bar{e}_{j})\in\mathfrak{k}_{\mathbb{C}}$. The operator
$D(u,\bar{u})$ has real spectrum for each $u\in V$. Hence, by
(\ref{LieInP}), $H_{j}\in i\mathfrak{k}$. Furthermore, the
subspace
\begin{equation}
\mathfrak{t}_{-}=i\bigoplus_{j=1}^{r}\mathbb{R}H_{j}\subset\mathfrak{k} \label{t-}%
\end{equation}
is abelian. Let $\mathfrak{t}=\mathfrak t_{-}\oplus\mathfrak t_{+}$ be a
Cartan subalgebra of $\mathfrak k$ and $\mathfrak{g}$ containing
$\mathfrak{t}_{-}$. Let
$\Delta(\mathfrak{g}_{\mathbb{C}},\mathfrak{t}_{\mathbb{C}})$ be
the set of roots of $\mathfrak{t}_{\mathbb{C}}$ in
$\mathfrak{g}_{\mathbb{C}}$. Recall that a root $\alpha\in\Delta$
is called \textit{compact} if $(\mathfrak{g}_{\mathbb{C}})_{\alpha}
\subset\mathfrak{k}_{\mathbb{C}}$ and \textit{non-compact} if $(\mathfrak{g}%
_{\mathbb{C}})_{\alpha}\subset\mathfrak{p}_{\mathbb{C}}$. Denote
by $\Delta_{c}$ and  $\Delta_{n}$ the set of compact,
respectively, non-compact roots. Finally, we recall that two roots
$\alpha\not =\pm\beta$ are called
\textit{strongly orthogonal} if $\alpha\pm\beta\not \in\Delta(\mathfrak{g}%
_{\mathbb{C}},\mathfrak{t}_{\mathbb{C}})$. We choose
$\gamma_{j}\in(\mathfrak t^{\ast})_{\mathbb{C}}$ so that
\[
\gamma_{j}(H_{k})=2\delta_{jk}%
\]
and  vanishes on $(\mathfrak t_{+})_{\mathbb{C}}$. Then
$\{\gamma_{1},\dots,\gamma_{r}\}$ is a maximal set of strongly
orthogonal non-compact roots. We order them so that
\[
\gamma_{1}>\dots>\gamma_{r}.
\]
The element $e:=e_{1}+\ldots+e_{r}$ is a \textit{maximal tripotent} and
\begin{equation}
Z_{0}=\frac{1}{2}D(e,\bar{e})=\frac{1}{2}\sum_{r=1}^{r}H_{j}, \label{Z0}%
\end{equation}
where $Z_{0}$ corresponds to the identity map as before. In this notation the
genus of $\mathcal{D}$ is given by
\begin{equation}
p=p(\mathcal{D})=\frac{2d}{r}\,. \label{eq-gen}%
\end{equation}

\subsection{Real bounded symmetric domain}

Let $e=e_{1}+\ldots+ e_{r}$ be the maximal tripotent as before.
The map $\tau:V\rightarrow V$, $\tau(v)=Q(e)\bar{v}$ is a
conjugate linear involution of $V$ and induces a decomposition
$V=J\oplus iJ$, where $J$ is a formally real Euclidean Jordan
algebra with identity $e$. Notice that $\tau$ also defines an
involution, also denoted by $\tau$, of $\mathfrak{g}$ and $G$. The
involution on $G$ is simply given by
$[\tau(g)](z)=\tau(g(\tau(z)))$, for $z\in \mathcal{D}$. As $\tau$
commutes with the Cartan involution $\theta$ we have the following
decomposition of $\mathfrak{g}$ into eigenspaces:
\begin{align*}
\mathfrak{g}  &  =\mathfrak{h}\oplus\mathfrak{q}\\
&  =\mathfrak{h}_{k}\oplus\mathfrak{q}_{k}\oplus\mathfrak{h}_{p} \oplus\mathfrak{q}_{p}\\
&  = \mathfrak{k}_{h} \oplus\mathfrak{k}_{q}\oplus\mathfrak{p}_{h}\oplus\mathfrak{p}_{q},%
\end{align*}
where $\mathfrak{h}$ is the $+1$-eigenspace and $\mathfrak{q}$ is
the $-1$-eigenspace of $\tau$. The index $k$ indicates
intersection with $\mathfrak{k}$, etc. The restriction of $\theta$
to $\mathfrak{h}$ defines a Cartan involution on $\mathfrak{h}$
and the corresponding Cartan decomposition is
$\mathfrak{h}=\mathfrak{k}_{h}\oplus\mathfrak{p}_{h}$. The space
$\mathfrak{p}_{q}$ is then a real subspace of $\mathfrak{p}$.
Recall the definition $\xi_{v}(z):=v-Q(z)\bar{v}$. Then
\[
\mathfrak{p}_{h}=\{\xi_{v}\mid v\in J\}.
\]
Let $\mathcal{D}_{\mathbb{R}}=\mathcal{D}\cap
J=\mathcal{D}^{\tau}$. Then $\mathcal{D}_{\mathbb{R}}$ is a real
bounded symmetric domain. For the structure of these domains see
\cite{HO97} and \cite{Loos-bsd}. Let
\[
G(\mathcal{D}_{\mathbb{R}})=\left\{  g\in G\mid g(\mathcal{D}_{\mathbb{R}%
})=\mathcal{D}_{\mathbb{R}}\right\}  \,.
\]
Then $G(\mathcal{D}_{\mathbb{R}})$ is a closed subgroup of $G$
with finitely many connected components. Let
$H=G(\mathcal{D}_{\mathbb{R}})_{o}$ be the connected component
containing the identity $1\in G$. Since $e\in
\mathcal{D}_{\mathbb{R}}$ the subgroup $\{k\in K\mid k\cdot e=e\}$
is a maximal compact subgroup of $G(\mathcal{D}_{\mathbb{R}})$
and equals
$G(\mathcal{D}_{\mathbb{R}})\cap K$.\ By replacing $G(\mathcal{D}_{\mathbb{R}%
})$ by $H$ it follows that $H\cap K$ is a maximal compact subgroup of $H$ and
\[
\mathcal{D}_{\mathbb{R}}=H/H\cap K=G(\mathcal{D}_{\mathbb{R}})/G(\mathcal{D}%
_{\mathbb{R}})\cap K\,
\]
is a Riemannian symmetric space with the Bergman metric restricted
to $\mathcal{D}_{\mathbb{R}}$. Moreover, it is a totally geodesic
submanifold. Finally, $H=G^{\tau}_{o}=\left\{  a\in
G\mid\tau(a)=a\right\}  _{o}$. Thus $(G,H)$ is a symmetric pair.
We note that the group $H$ is not semisimple since
$\operatorname{exp}(\mathbb{R}\xi_{e})$ is in the center of $H$.
The group $H$ is invariant under the Cartan involution $\theta$
and $\theta|_{H}$ is a Cartan involution on $H$.

\subsection{Cayley transform}

In this subsection we discuss the realization of $\mathcal{D}$ as
a tube type domain $T(\Omega)=iJ+\Omega$. Let $\Omega$ be the cone
of positive elements in $J$: \
\[
\Omega=\{x^{2}\mid x\in J,\,\operatorname{det}(x)\neq0\}.
\]
Then $\Omega$ is a symmetric convex cone. Let
\[
T(\Omega)=iJ+\Omega=\{w\in V\mid\mathrm{Re}(w)\in\Omega\}\,.
\]
If $e-z$ is invertible in $V$ let
\begin{equation}
\gamma(z)=\frac{e+z}{e-z}=(e+z)(e-z)^{-1} \label{cayley1}%
\end{equation}
be the Cayley transform. Its inverse is
\begin{equation}
\gamma^{-1}(w)=(w-e)(w+e)^{-1} \label{cayley2}.
\end{equation}
The following is well known:

\begin{lemma}
Let the notation be as above. Then the following hold:
\begin{enumerate}
\item  The Cayley transform $\gamma$ is a biholomorphic transformation from
$\mathcal{D}$ into the \textit{Siegel domain} $T(\Omega)$.

\item  The Cayley transform $\gamma:\mathcal{D}\mapsto T(\Omega)$ maps the
real form $\mathcal{D}_{\mathbb{R}}$ of $\mathcal{D}$ onto the real form
$\Omega$ of $T(\Omega)$.

\item  The following diagram commutes
\[%
\begin{array}
[c]{ccc}%
\mathcal{D}_{\mathbb{R}} & \overset{\gamma}{\longrightarrow} & \Omega\\
\downarrow &  & \downarrow\\
\mathcal{D} & \overset{\gamma}{\longrightarrow} & T(\Omega)
\end{array}
\]
where the vertical arrows are the inclusion maps.
\end{enumerate}
\end{lemma}

Notice that the Cayley transform can be realized by an element,
also denoted  $\gamma$, in $G_{\mathbb{C}}$. In fact, there exists
a $X\in\mathfrak{q}_{p}$ such that $\mathrm{ad}(X)$ has
eigenvalues $0,1,$ and $-1$ and $\gamma
=\operatorname{exp}(\frac{\pi i}{4}X)$. Define $G^{\gamma}:=\gamma
G\gamma^{-1}$ and $H^{\gamma}=(K_{\mathbb{C}}\cap
G^{\gamma})_{o}$. Then
\begin{equation}
G^{\gamma}(\Omega)_{o}=\left\{  g\in G^{\gamma}\mid g(\Omega)=\Omega\right\}
=H^{\gamma}=\gamma H \gamma^{-1}. \label{eq-Hgamma}%
\end{equation}
Let $\mathfrak{g}^{\gamma}:=\mathrm{Ad}(\gamma)(\mathfrak{g})$ be
the Lie algebra of $G^{\gamma}$. We collect some important facts
in the following lemmas:

\begin{lemma}
\label{le-ggamma}Let the notation be as above. Then the following
holds:
\begin{equation}
\mathfrak{g}^{\gamma}=\mathfrak{h}_{k}+i\mathfrak{h}_{p}+i\mathfrak{q}_{k}+\mathfrak{q}_{p}
\label{eq-ggamma}%
\end{equation}
\end{lemma}
Thus $(\mathfrak{g}^{\gamma},\mathfrak{k}^{\gamma}=\mathfrak{h}_{k}+i\mathfrak{h}_{p}%
,\mathfrak{h}^{\gamma}=\mathfrak{h}_{k}+i\mathfrak{q}_{k})$ is the
\textit{Riemannian dual} of
$(\mathfrak{g},\mathfrak{k},\mathfrak{h})$.

\begin{lemma}
With notation as above we have
$$\mathfrak{h}^{\gamma}=\mathfrak{k}_{\mathbb{C}}\cap\mathfrak{g}^{\gamma} \text{ and }
 Z_{0} \text{ is central in } \mathfrak{h}^{\gamma}.$$
\end{lemma}

We refer to \cite{HO97} for further information. The group $G^{\gamma}$ acts
on $T(\Omega)$ by transforming the action of $G$:%
\[
g\cdot z=\gamma((\gamma^{-1}g\gamma)\cdot(\gamma^{-1}(z))).
\]
For simplicity we will sometimes write $\gamma$ for the adjoint
action of $\gamma$ on $\mathfrak{g}_{\mathbb{C}}$ and the
conjugation with $\gamma$ in $G_{\mathbb{C}}$.

\subsection{Roots}

In this subsection we describe the structure of restricted roots
for the group $H$. Let $\xi_{j}=\xi_{e_{j}}$ and
\[
\xi=\xi_{1}+\dots+\xi_{r}%
\]
By $\mathrm{SU}(1,1)$-reduction we have the following:

\begin{lemma}
\label{le-xi}$\xi_{j}=\mathrm{Ad}(\gamma)^{-1}(H_{j})$ and $\xi=2\mathrm{Ad}%
(\gamma)^{-1}(Z_{0})=-\mathrm{Ad}(\gamma)(2Z_{0})$.
\end{lemma}

Let
\[
\mathfrak{a}=\bigoplus_{j=1}^{r}\mathbb{R}\xi_{j}\,,
\]
and define $\beta_{j}\in\mathfrak{a}^{\ast}$ by
$\beta_{j}(\xi_{i})=2\delta_{ij}$. Then $\mathfrak{a}$ is maximal
abelian in $\mathfrak{h}_{p}$. In fact, $\mathfrak{a}$ is also
maximal abelian in $\mathfrak{p}$. We remark that $\gamma^2
H_j=-H_j$, $\gamma(\mathfrak t_{-})=\gamma^{-1}(\mathfrak
t_{-})=\mathfrak a$ and
$\beta_{j}=\gamma_{j}\circ\mathrm{Ad}(\gamma)$.

We will often identify $\mathfrak{a}_{\mathbb{C}}^{\ast}$ with
$\mathbb{C}^{r}$ using the map
$\mathbf{\alpha\mapsto}{\mathbf{\alpha}}=\alpha_{1}\beta
_{1}+\ldots+\alpha_{r}\beta_{r}$. The root system
$\Delta(\mathfrak h,\mathfrak a)$ is of type A:
\[
\Delta(\mathfrak{h},\mathfrak{a})=\pm\left\{
\frac{\beta_{j}-\beta_{k}}{2}\mid1\leq j\neq k\leq r\right\}  .
\]
We fix an ordering of the roots so that
\[
\beta_{1}>\beta_{2}>\dots>\beta_{r}.
\]
Then the corresponding system $\Delta^{+}=\Delta^{+}(\mathfrak h,\mathfrak a)$
of positive roots is given by
\[
\Delta^{+}=\left\{  \frac{\beta_{j}-\beta_{k}}{2}\mid1\leq j<k\leq r\right\}
.
\]
The root spaces $\mathfrak{h}_{(\beta_{j}-\beta_{k})/2}$ all have
the same
dimension which we denote by $a$. Then the half sum of the positive roots is%

\begin{equation}
\rho=\sum_{j=1}^{r}\rho_{j}\beta_{j}=\frac{a}{4}\sum_{j=1}^{r}((r+1)-2j)\beta
_{j} \label{rho}.%
\end{equation}

\begin{example}
\label{ex-su11}In the case of $\mathcal{D}=\left\{
z\in\mathbb{C}\mid\left| z\right|  <1\right\}  $ and
$G=\mathrm{SU}(1,1)$ acting on $\mathcal{D}$ in the usual way
\[
\left(
\begin{array}
[c]{cc}%
\alpha & \beta\\
\bar{\beta} & \bar{\alpha}%
\end{array}
\right)  \cdot z=\frac{\alpha z+\beta}{\bar{\beta}z+\bar{\alpha}}%
\]
our notation is:%
\begin{align*}
2Z_{0} &  =H_{1}=\left(
\begin{array}
[c]{cc}%
1 & 0\\
0 & -1
\end{array}
\right)  \,,\\
\gamma &  =\frac{1}{\sqrt{2}}\left(
\begin{array}
[c]{cc}%
1 & 1\\
-1 & 1
\end{array}
\right)  =\operatorname{exp}\left(  \frac{\pi i}{4}\left(
\begin{array}
[c]{cc}%
0 & -i\\
i & 0
\end{array}
\right)  \right)  ,\\
\xi &  =\left(
\begin{array}
[c]{cc}%
0 & 1\\
1 & 0
\end{array}
\right)  ,\\
\mathcal{D}_{\mathbb{R}} &  =(-1,1)\,,\\
J &  =\mathbb{R},\\
\Omega&=\mathbb{R}^{+}\mathbb{\,},\\
T(\Omega) &  =\left\{  z\in\mathbb{C}\mid\mathrm{Re}(z)>0\right\}  \,.\qquad
\end{align*}
\end{example}

\subsection{Conical functions, spherical functions, and invariant polynomials}

Let $\{e_{j}\}_{j=1}^{r}$ be the fixed frame as before. Let $u_{j}:=\sum
_{k=1}^{j}e_{k}$, $j=1,\cdots,r$. Let $V_{j}:=\{z\in V\mid D(u_{j},\bar{u}%
_{j})z=2z\}$. Then $V_{j}$ is a Jordan $\star$-subalgebra of $V$
with a determinant polynomial denoted by $\Delta_{j}$. We extend
$\Delta_{j}$ to all of $V$ via
$\Delta_{j}(z):=\Delta_{j}(\mathrm{pr}_{V_{j}}(z))$, where
$\mathrm{pr}_{V_{j}}$ is the orthogonal projection onto $V_{j}$.
The polynomials $\Delta_{j}$ are called (principal)
\textit{minors}. Notice that
$\Delta_{r}(w)=\Delta(w)=\operatorname{det}(w)$. For any
$\mathbf{\alpha}=(\alpha_{1},\cdots,\alpha_{r})\in{\mathbb{C}}^{r}$
consider the associated \textit{conical function} (\cite{FK-book},
p. 122):
\[
\Delta_{\mathbf{\alpha}}(w):=\Delta_{1}^{\alpha_{1}-\alpha_{2}}(w)\Delta
_{2}^{\alpha_{2}-\alpha_{3}}(w)\cdots\Delta_{r}^{\alpha_{r}}(w),\,\,w\in V.
\]
Notice that if \ $w=\sum_{j=1}^{r}w_{j}e_{j}$ \ then \ $\Delta_{\alpha
}(w)=\prod_{j=1}^{r}w_{j}^{\alpha_{j}}$. Thus the conical functions are
generalizations of the power functions.

Let  $L:=H\cap K\subset G\cap G^{\gamma}$ and define
\begin{equation}
\psi_{\mathbf{\alpha}}(z):=\int_{L}\Delta_{\mathbf{\alpha}}(lz)dl\,,\qquad
z\in V \label{psi}%
\end{equation}
The function $\psi_{\lambda+\rho}$ is the spherical function on
$\Omega$ corresponding to $\lambda$ (\cite{FK-book}, Theorem XIV.
3.1).  We identify $(\mathfrak{t}_{-}^{\ast})_{\mathbb{C}}$ with
$\mathbb{C}^{r}$ via
$m_{1}\gamma_{1}+\ldots+m_{r}\gamma_{r}\leftrightarrow(m_{1},\ldots,m_{r}%
)\in\mathbb{C}^{r}.$ If $\lambda={\mathbf{m}}$, where $\mathbf{m}=(m_{1}%
,\dots,m_{r})$ is a tuple of non-negative integers such that
$m_{1}\geq m_{2}\geq\dots\geq m_{r}\geq0$, then the functions
$\Delta_{\lambda}$ and $\psi_{\lambda}$ are holomorphic
polynomials on the whole space $V$ and $\psi_{\lambda}$ is
$L$-invariant. Let $\mathcal{P}(V)$ be the space of holomorphic
polynomials on $V$, considered as $K$-space by the regular action,
$k\cdot p(z)=p(k^{-1}\cdot z)$. Let
\begin{equation}
\mathbf{\Lambda}=\{\mathbf{m}\in\mathbb{N}_{0}^{r}\mid m_{1}\geq m_{2}%
\geq\ldots\geq m_{r}\geq0\}\,. \label{eq-lambda}%
\end{equation}
Then we have the well known Schmid decomposition (\cite{Schmid}
and \cite{FK-book}, Theorem XI.2.4).

\begin{lemma}
The space of polynomials $\mathcal{P}(V)$ decomposes as a
$K$-representation into
\begin{equation}
\mathcal{P}(V)=\sum_{\mathbf{m}\in\mathbf{\Lambda}}\mathcal{P}_{\mathbf{m}%
},\label{schmid}%
\end{equation}
where each $\mathcal{P}_{\mathbf{m}}$ is of lowest weight $-\mathbf{m}%
=-(m_{1}\gamma_{1}+\dots+m_{r}\gamma_{r})$, with $m_{1}\geq\cdots\geq
m_{r}\geq0$ being integers. In this case the polynomial $\Delta_{\mathbf{m}}$
is a lowest weight vector in $\mathcal{P}_{\mathbf{m}}$ and $\psi
_{{\mathbf{m}}}$ is up to constants the unique $L$-invariant polynomials in
$\mathcal{P}_{\mathbf{m}}$. In particular
\begin{equation}
\mathcal{P}(V)^{L}=\sum_{\mathbf{m}\in\mathbf{\Lambda}}\mathbb{C}%
\psi_{{\mathbf{m}}}.\label{schmid-L}%
\end{equation}
\end{lemma}

As $\psi_{\lambda+\rho}$ is spherical on $\Omega$ and because the
Cayley transform commutes with the $L$-action we get the following
lemma.

\begin{lemma}
\label{l2.2} Let $\lambda\in\mathfrak{a}_{\mathbb{C}}^{\ast}$.
Then the spherical function $\phi_{\lambda}(x)$ on the real
bounded symmetric domain $\mathcal{D}_{\mathbb{R}}$ is given by
\[
\phi_{\lambda}(x)=\psi_{{i\lambda}+\rho}\left(
\frac{e+x}{e-x}\right) =\psi_{i\lambda+\rho}\circ\gamma(x)\,,\quad
x\in\mathcal{D}_{\mathbb{R}}\,.
\]
\end{lemma}

\subsection{The $\Gamma$-function on symmetric cones}

The $H$-invariant measure on $\Omega$ is given by
\[
d\mu_0=\Delta(x)^{-\frac{d}{r}}dx,\label{inv-measure}%
\]
where $d=\dim_{\mathbb{R}}(J)=\dim_{\mathbb{C}}(V)$. The
\textit{Gindikin-Koecher Gamma function} associated with the
convex, symmetric cone $\Omega$ is defined by
\[
\Gamma_{\Omega}(\lambda):=\int_{\Omega}e^{-\mathrm{Tr}(x)}\Delta_{\lambda
}(x)\Delta(x)^{-d/r}dx,
\]
The integral converges if and only if
$\mathrm{Re}(\lambda_{j})>(j-1)a/2$ for \ $j=1,2,\cdots,r$,
(\cite{FK-book}, Theorem VII.1.1.). Moreover, the convergence is
absolute and uniform on compact subsets of $\mathbb{C}^n$. Using
the identification  $\mathfrak{a}_{\mathbb{C}}^{\ast}$ with
$\mathbb{C}^n$ via the map
$(\lambda_1,\ldots,\lambda_n)\leftrightarrow
\lambda_1\beta_1+\cdots \lambda_n\beta_n$ we have
\[
\Gamma_{\Omega}(\lambda)=\left(  2\pi\right)  ^{\frac{d-r}{2}}\prod_{j=1}%
^{r}\Gamma\left(  \lambda_{j}-(j-1)\frac{a}{2}\right),
\]
where $\Gamma$ is the usual Gamma function. In particular, it
follows that $\Gamma_{\Omega}$ has a meromorphic continuation to all of $\mathfrak{a}%
_{\mathbb{C}}^{\ast}$. We will view $\Gamma_{\Omega}$ as a
function on $\mathbb{C}^{r}$ and on
$\mathfrak{a}_{\mathbb{C}}^{\ast}$ using our identification
above. We also adopt the
notation $\beta_{0}=\sum_{j=1}^{r}\beta_{j}$ and
$\Gamma_{\Omega}(\nu )=\Gamma_{\Omega}(\nu\beta_{0})$, where $\nu
\in \mathbb{C}$. Finally, we define
\[
\left(  \lambda\right)  _{\mathbf{m}}=\frac{\Gamma_{\Omega}(\lambda+\mathbf{m}%
)}{\Gamma_{\Omega}(\lambda)}\,.
\]

\subsection{The Laplace transform}

We recall here a few facts about the Laplace transform on
$\Omega$. Let $\mu$ be a (complex) Radon measure on $\Omega$ such
that $x\mapsto e^{-(t\mid x)}\in L^{1}(\Omega,d|\mu|)$ for all
$t\in\Omega$. Define the \textit{Laplace transform of} $\mu$ by
\[
\mathcal{L}(\mu)(w):=\int_{\Omega}e^{-(w\mid
x)}\,d\mu(x)=\int_{\Omega }e^{-(t\mid x)}e^{-i(s\mid x)}d\mu(x),
\]
for $w=t+is\in T(\Omega)$. Then $\mathcal{L}(\mu)$ is holomorphic
on $T(\Omega)$. In particular, if $f\in L^{1}(\Omega,d\mu)$ then
$f(x)d\mu$ is a finite measure and hence
$\mathcal{L}(f):=\mathcal{L}(fd\mu)$ is well defined. Furthermore,
if $\nu\in\mathbb{C}$ is such that $\mathrm{Re}(\nu)>(r-1)\frac
{a}{2}$ let $d\mu_{\nu}(x)=$ $\Delta(x)^{\nu-d/r}\,dx$. Notice
that $\Delta_{\nu-d/r}(x)=\Delta(x)^{\nu-d/r}$. Then $d\mu_{\nu}$
is a quasi-invariant measure on $\Omega$ and the $H$-invariant
measure corresponds to $d\mu_{\circ}$. If $\nu>0$ is real, then we
let
\begin{equation}
L_{\nu}^{2}(\Omega)=L^{2}(\Omega,d\mu_{\nu})\,. \label{L2nu}%
\end{equation}
We define $\mathcal{L}_{\nu}:L_{\nu}^{2}(\Omega)\rightarrow\mathcal{O}%
(T(\Omega))$ by
\[
\mathcal{L}_{\nu}(f)(z)=\mathcal{L}(fd\mu_{\nu})(z)=\int_{\Omega}e^{-(z\mid
x)}f(x)\,d\mu_{\nu}(x), \quad, z\in T(\Ome)
\]
Here $\mathcal{O}(T(\Omega)$ denotes the space of holomorphic
functions on $T(\Omega)$. We notice the following which follows
directly from Proposition VII.1.2 of \cite{FK-book}, p. 124 by the
holomorphicity of both sides in $w$:

\begin{lemma}
Let $\nu\in\mathbb{C}$ be such that $\mathrm{Re}(\nu)>(r-1)\frac{a}{2}$. Then
for any $w\in T(\Omega)$%
\[
\mathcal{L}(\mu_{\nu})(w)=\Gamma_{\Omega}(\nu)\Delta(w)^{-\nu}.
\]
\end{lemma}

\subsection{Unitary highest weight representations}

In this subsection we review some simple facts on scalar valued
unitary highest weight representations. We restrict the discussion
to what we will need later on. From now on $v\mapsto \bar{v}$ denotes
conjugation with respect to the real form $J$. Let
$\tilde{G}^{\gamma}$ be the universal covering group of
$G^{\gamma}$. Then $\tilde{G}^{\gamma}$ acts on $T(\Omega)$ by
$(g,z)\mapsto\kappa(g)\cdot z$ where
$\kappa:\tilde{G}^{\gamma}\rightarrow G^{\gamma}$ is the canonical
projection. For $\nu>1+a(r-1)$ let $\mathcal{H}_{\nu}(T(\Omega))$
be the space of holomorphic functions
$F:T(\Omega)\rightarrow\mathbb{C}$ such that
\begin{equation}
\left|  \left|  F\right|  \right|  _{\nu}^{2}:=\alpha_{\nu}\int_{T(\Omega
)}|F(x+iy)|^{2}\Delta(y)^{\nu-2d/r}\,dxdy<\infty\label{eq-normub}%
\end{equation}
where
\begin{equation}
\alpha_{\nu}=\frac{2^{r\nu}}{(4\pi)^{d}\Gamma_{\Omega}(\nu-d/r)}\,.
\label{eq-cnu}%
\end{equation}
Then $\mathcal{H}_{\nu}(T(\Omega))$ is a non-trivial Hilbert
space. For $\nu\leq1+a(r-1)$ this space reduces to $\{0\}$. If
$\nu=2d/r$ this is the \textit{Bergman space}. For
$g\in\tilde{G}^{\gamma}$ and $z\in T(\Omega)$, let
$J_{g}(z)=J(g,z)$ be the \textit{complex} Jacobian determinant of
the action of $\tilde{G}^{\gamma}$ on $T(\Omega)$ at the point
$z$. We will use the same notation for elements $g\in G$ and
$z\in\mathcal{D}$. Then
\[
J(ab,z)=J(a,b\cdot z)J(b,z)
\]
for all $a,b\in\tilde{G}^{\gamma}$ and $z\in T(\Omega)$. It is well known that
for $\nu>1+a(r-1)$ that
\begin{equation}
\pi_{\nu}(g)f(z)=J(g^{-1},z)^{\nu/p}f(g^{-1}\cdot z) \label{eq-actionub}
\end{equation}
defines a unitary irreducible representation of
$\tilde{G}^{\gamma}$. In \cite{Rossi-Vergne},
\cite{Wallach} and \cite{FK} it was shown that this unitary
representation $(\pi_{\nu },\mathcal{H}_{\nu}(T(\Omega)))$ has an
analytic continuation to the half-interval $\nu>(r-1)\frac{a}{2}$.
Here the representation $\pi_{\nu}$ is given by the same formula
(\ref{eq-actionub}) but the formula for the \textit{norm} in
(\ref{eq-normub}) is no longer  valid. We collect the necessary
information from \cite{FK-book}, p. 260, in particular, Theorem
XIII.1.1. and Proposition XIII.1.2., in the following theorem. We
will give a new proof of (4) later using only part (3).

\begin{theorem}
\label{th-2.5}Let the notation be as above. Assume that for
$\nu>1+a(r-1)$ then the following hold:

\begin{enumerate}
\item  The space $\mathcal{H}_{\nu}(T(\Omega))$ is a reproducing Hilbert space.

\item  The reproducing kernel of $\mathcal{H}_{\nu}(T(\Omega))$ is given by
\[
K_{\nu}(z,w)=\Gamma_{\Omega}(\nu)\Delta\left(  z+\bar{w}\right)  ^{-\nu}%
\]

\item  If $\nu>(r-1)\frac{a}{2}$ then there exists a Hilbert space
$\mathcal{H}_{\nu}(T(\Omega))$ of holomorphic functions on
$T(\Omega)$ \ such that $K_{\nu}(z,w)$ defined in (2) is the
reproducing kernel of that Hilbert
space. The group $\tilde{G}^{\gamma}$ acts unitarily on $\mathcal{H}_{\nu}%
(T(\Omega))$ by the action defined in (\ref{eq-actionub}).

\item  The map
\[
L_{\nu}^{2}(\Omega)\ni f\mapsto F=\mathcal{L}_{\nu}(f)\in\mathcal{H}_{\nu
}(T(\Omega))
\]
is a unitary isomorphism and

\item  If $\nu>(r-1)\frac{a}{2}$ then the functions
\[
q_{\mathbf{m},\nu}(z):=\Delta(z+e)^{-\nu}\psi_{\mathbf{m}}\left(  \frac
{z-e}{z+e}\right),  \,\qquad\mathbf{m}\in\mathbf{\Lambda},%
\]
form an orthogonal basis of $\mathcal{H}_{\nu}(T(\Omega))^{L}$,
the space of $L$-invariant functions in
$\mathcal{H}_{\nu}(T(\Omega))$.
\end{enumerate}
\end{theorem}

\subsection{Weighted Bergman spaces on the bounded symmetric domain}

In the last section we concentrated on the unbounded realization. We will now
shift our attention to the bounded realization $\mathcal{D}$. Let
\begin{equation}
h(z,w):=\Delta(e-z\bar{w}).
\end{equation}
We note the $h(z, w)^{-p}$ is  the \textit{Bergman kernel} on
$\mathcal{D}$. We collect a few well known facts about $h(z,w)$ in
the following lemma.
\begin{lemma}
\label{le-proph} Let $g\in G$ and $z,w\in\mathcal{D}$. Then the
following hold:
\begin{enumerate}
\item $h(z,w)$ is holomorphic in the first variable and anti-holomorphic in
the second variable.
\item $h(w,z)=\overline{h(z,w)}$ for all $z,w\in\mathcal{D}$.
\item $h(g\cdot z,g\cdot w)=J(g,z)^{1/p}\overline{J(g,w)^{1/p}}\,h(z,w)$. In
particular, $h(z,z)=|J(g,0)|^{2/p}>0$ where $g\in G$ is chosen
such that $g\cdot0=z$.
\item  If $x,y\in\mathcal{D}_{\mathbb{R}}$ then $h(x,y)>0$. In particular, all
powers $h(x,y)^{\mu}$, $\mu\in\mathbb{R}$,
$x,y\in\mathcal{D}_{\mathbb{R}}$, are well defined.
\end{enumerate}
\end{lemma}

Let $h(z)=h(z,z)=\Delta(e-z\bar{z})$. Then the measure
\begin{equation}
dm_{\nu}(z)=h(z)^{\nu-p}dz
\end{equation}
is a quasi-invariant measure on $\mathcal{D}$ and $dm_{0}(z)$ is
invariant. Furthermore,
\begin{equation}
d\eta(x)=h(x)^{-p/2}dx
\end{equation}
is an $H$-invariant measure on $\mathcal{D}_{\mathbb{R}}$. Let $\tilde{G}$ be
the universal covering group for $G$. As for $T(\Omega)$ there exists a
Hilbert space $\mathcal{H}_{\nu}(\mathcal{D})$ of holomorphic functions
$\mathcal{D}$ with a unitary $\tilde{G}$-action, given by
\[
\pi_{\nu}(g)f(z)=J(g^{-1},z)^{\frac{\nu}{p}}f(g^{-1}z)
\]
if $\nu>(r-1)a/2$. Here $z\mapsto J(g,z)$ is again the Jacobian of the action
of $G$ on $V$. If $\nu>1+a(r-1)$ then the norm on $\mathcal{H}_{\nu
}(\mathcal{D})$ is given by
\begin{equation}
\Vert F\Vert_{\nu}^{2}=\Vert
F\Vert_{\mathcal{H}_\nu(\mathcal{D})}^{2}:=d_{\nu}\int_{\mathcal{D}}|F(z)|^{2}\,dm_{\nu}(z)
\end{equation}
where
\[
d_{\nu}=\frac{1}{\pi^{d}}\frac{\Gamma_{\Omega}(\nu)}{\Gamma_{\Omega}(\nu
-d/r)}\,.
\]
The constant $d_{\nu}$ is chosen so that the constant function $z\mapsto1$ has
norm one.

\begin{lemma}
\label{le-HnuD} Let the notation be as above. Then the following
hold:
\begin{enumerate}
\item  If $F\in\mathcal{H}_{\nu}(T(\Omega))$ then the function
\[
\pi_\nu(\gamma^{-1})(F)(w)=2^\frac{r\nu}{2}
\Delta(e-w)^{-\nu}F\circ\gamma(w)=\Delta(e-w)^{-\nu
}F((e+w)(e-w)^{-1})
\]
belongs to $\mathcal{H}_{\nu}(\mathcal{D})$ and
\[
\pi_\nu(\gamma^{-1}):\mathcal{H}_{\nu}(T_{\Omega})\rightarrow\mathcal{H}_{\nu
}(\mathcal{D})
\]
is a linear isomorphism onto $\mathcal{H}_{\nu}(\mathcal{D})$.

\item  The inverse $\pi_\nu(\gamma):\mathcal{H}_{\nu}(\mathcal{D}%
)\rightarrow\mathcal{H}_{\nu}(T(\Omega))$ is given by
\[
\pi_\nu(\gamma)(F)(z)=2^\frac{r \nu }{2}\Delta\left(  z+e\right)  ^{-\nu}%
F\circ\gamma^{-1}(z)=2^\frac{r \nu }{2}\Delta\left(  z+e\right)
^{-\nu}F(\left( z-e)(z+e)^{-1}\right)  \,.
\]

\item  Let $F\in\mathcal{H}_{\nu}(T(\Omega))$ then
\[
\Vert\pi_\nu(\gamma^{-1})(F)\Vert_{\mathcal{H}_{\nu}(\mathcal{D})}^{2}=
\Gamma_{\Omega}(\nu)\Vert F\Vert_{\mathcal{H}_{\nu}(T(\Omega
))}^{2}%
\]

\item  If $g\in G$ and $F\in\mathcal{H}_{\nu}(\mathcal{D})$ then
\[
\pi_\nu(\gamma)(\pi_{\nu}(g)F)=\pi_{\nu}(\gamma
g\gamma^{-1})\pi_\nu(\gamma)(F)\,.
\]

\item  If $\nu>(r-1)\frac{a}{2}$ then $P(V)\subset\mathcal{H}_{\nu
}(\mathcal{D})$ and $P(V)$ is dense in $\mathcal{H}_{\nu}(\mathcal{D})$.

\item  If $\nu>(r-1)\frac{a}{2}$ then the functions $\psi_{\mathbf{m}}$,
$\mathbf{m}\in\mathbf{\Lambda}$, form an orthogonal basis for $\mathcal{H}%
_{\nu}(\mathcal{D})^{L}$.
\end{enumerate}
\end{lemma}

\begin{proof}
That the map is an isomorphism follows from \cite{FK-book},
Proposition XIII,1.3. The intertwining relation in property (4) is
a simple calculation and is left to the reader.
\end{proof}

In particular, it follows that (\cite{FK-book}, Proposition
XIII.1.4):

\begin{lemma}
The reproducing kernel of $\mathcal{H}_{\nu}(\mathcal{D})$ is given by
$K_{\nu}(z,w)=h(z,w)^{-\nu}=\Delta(e-z\bar{w})^{-\nu}$.
\end{lemma}

Define the \textit{Fock-space} on $V$ to be the space of
holomorphic functions on $V$ such that
\begin{equation}
\left|  \left|  F\right|  \right|  ^{2}_{\mathcal{F}(V)}=\pi^{-d}\int
_{V}|F(z)|^{2}\,e^{-\|z\|^{2}}\,dz<\infty\,.
\end{equation}
The following result is proved by Faraut and Koranyi \cite{FK}, see also
\cite{FK-book}, Proposition XI.4.1 and Proposition XIII.2.2. This will play an
essential role in our work.

\begin{theorem}
Assume that $\nu>1+a(r-1)$ and $\mathbf{m}\geq0$. Then the norms
of $\psi_{\mathbf{m}}$ in the Fock-space and weighted Bergman
spaces are given by
\[
\Vert\psi_{\mathbf{m}}\Vert_{\mathcal{F}(V)}^{2}=\frac{1}{d_{\mathbf{m}}%
}\left(  \frac{d}{r}\right)  _{\mathbf{m}}%
\]
and
\[
\Vert\psi_{\mathbf{m}}\Vert_{\mathcal{H}_{\nu}(\mathcal{D})}^{2}=\frac
{1}{d_{\mathbf{m}}}\frac{\left(  \frac{d}{r}\right)  _{\mathbf{m}}}%
{(\nu)_{\mathbf{m}}}%
\]
respectively, where $d_{\mathbf{m}}$ is the dimension of the space
$\mathcal{P}_{\mathbf{m}}$.
\end{theorem}

\begin{corollary}
\label{cor-orthonormal} Assume that $\nu>1+a(r-1)$. Then the functions
\[
\sqrt{\frac{d_{\mathbf{m}}(\nu)_{\mathbf{m}}}{\left(  \frac{d}{r}\right)
_{\mathbf{m}}}}\,\psi_{\mathbf{m}}\,,\quad\mathbf{m}\in\mathbf{\Lambda}%
\]
form an orthonormal basis for $\mathcal{H}_{\nu}(\mathcal{D})^L$.
\end{corollary}

\begin{proof}
This follows from Lemma \ref{le-HnuD}, part 6, and the above
Theorem.
\end{proof}

\section{Berezin transform and generalized Segal-Bargmann transform for
$\mathcal{D}_{\mathbb{R}}$}

\subsection{Restriction to $D$ of holomorphic functions on $\mathcal{D}$}

In this subsection we discuss the restriction principle for the
bounded symmetric space $\mathcal{D}_{\mathbb{R}}$,
(\cite{o00,OOe96,gkz-bere-rbsd}) .  In particular, we give an
exact bound for the parameter $\nu$ such that the restriction map
$R$, defined below,  maps $\mathcal{P}(V)$ into $L^{2}(H/L)$.

Recall that the $H$-invariant measure on
$\mathcal{D}_{\mathbb{R}}$ is given by
$d\eta(x)=h{(x)^{-\frac{p}{2}}}dx$. The group $H$ acts unitarily
on $L^{2}(\mathcal{D}_{\mathbb{R}},d\eta)$ by
$g\cdot f(x)=f(g^{-1}\cdot x)$. Furthermore,
\[
L^{2}(\mathcal{D}_{\mathbb{R}},d\eta)\simeq_{H}
L^{2}(H/L,d\dot{h})
\]
where $d\dot{h}$ denotes an $H$-invariant measure on $H/L$.

\begin{lemma}
\label{le-sqrint} Let $\nu\in\mathbb{R}$. Then
$p(x)h(x)^{\nu/2}\in L^{2}(\mathcal{D}_{\mathbb{R}},d\eta)$, for
all $p\in\mathcal{P}(V)$, if and only if $\nu>\frac{a(r-1)}{2}$.
\end{lemma}

\begin{proof}
Let $p\in\mathcal{P}(V)$. As the closure of
$\mathcal{D}_{\mathbb{R}}$ is compact it follows that $p$ is
bounded on $\mathcal{D}_{\mathbb{R}}$. Hence it is enough to show
that $h^{\nu/2}\in L^{2}(\mathcal{D}_{\mathbb{R}},d\eta)$ if and
only if $\nu>(r-1)\frac {a}{2}$. Note that $h^{\nu/2}$ is
$L$-invariant. Writing the invariant measure
on $H/L$ using polar-coordinates gives for every $L$-invariant function:%

\[
\int_{H/L}f(h\cdot0)d\dot{h}=\int_{t_{1}>t_{2}>\ldots>t_{r}}%
f(\operatorname{exp}(\sum_{j=1}^{r}t_{j}\xi_{j})\cdot0)\left(  \prod_{i<j}%
\sinh(t_{i}-t_{j})\right)  ^{a}\,dt_{1}\ldots dt_{r}\,.
\]

If $G=\mathrm{SU}(1,1)$ then $\xi_{1}=\left(
\begin{array}
[c]{cc}%
0 & 1\\
1 & 0
\end{array}
\right)  $. Then
\[
g_{t}:=\operatorname{exp}(t\xi_{1})=\left(
\begin{array}
[c]{cc}%
\cosh(t) & \sinh(t)\\
-\sinh(t) & \cosh(t)
\end{array}
\right)
\]
and hence
\[
g_{t}\cdot0=\operatorname{tanh}(t)\,
\]
A similar calculation in the general case gives:
\[
\operatorname{exp}(\sum_{j=1}^{r}t_{j}\xi_{j})\cdot0=\sum_{j=1}^{r}%
\operatorname{tanh}(t_{j})e_{j}\,.
\]
Thus
\begin{align*}
h(\operatorname{exp}(\sum_{j=1}^{r}t_{j}\xi_{j})\cdot0) &  =\Delta
(e-\sum_{j=1}^{r}\operatorname{tanh}^{2}(t_{j})e_{j})\\
&  =\prod_{j=1}^{r}(1-\operatorname{tanh}^{2}(t_{j}))\\
&  =\prod_{j=1}^{r}\frac{1}{\cosh^{2}(t_{j})}\,.
\end{align*}
It follows that
\[
\int_{\mathcal{D}_{\mathbb{R}}}|h(x)^{\nu/2}|^{2}h(x)^{-p/2}\,dx=\int
_{t_{1}>\ldots>t_{r}}\prod_{j=1}^{r}\cosh(t_{j})^{-2\nu}\left(
\prod _{i<j}\sinh(t_{i}-t_{j})\right)  ^{a}\,dt_{1}\dots dt_{r}.
\]

Let $\phi (t)=\frac {1-e^{-2t}}{2}$ and observe that $\phi$ is
non-negative and increasing, $\phi(t)=0$ if and only if $t=0$,
$\phi$ is bounded (by $\frac 1 2$) on $[0,\infty)$, and $\sinh
(t_i-t_j)=e^{t_i-t_j}\phi(t_i-t_j)$. In a similar way let
$\psi(t)=\frac {1+e^{-2t}}{2}$. Observe that $\psi(t)\in(\frac 1 2
,1]$, for all $t\in [0,\infty)$,  and $\cosh
(t_j)=e^{t_j}\psi(t_j)$. Let $t=(t_1,\ldots,t_r)$ and define
$\Phi(t)=\Pi_{i<j}\phi(t_i-t_j)$ and $\Psi(t)=\Pi_{j=1}^r
\psi(t_j).$ In this notation we have
$$\int_{\mathcal{D}_{\mathbb{R}}}|h(x)^{\nu/2}|^{2}h(x)^{-p/2}\,dx=\int
_{t_{1}>\ldots>t_{r}}\Phi(t)\Psi(t)e^{\sum_{j-1}^r (-2\nu
+(r+1-2j)a)t_j} \,dt_{1}\dots dt_{r}.$$

Now suppose $\nu>(r-1)\frac a 2$. Then $-2\nu+(r+1-2j)a<0$ for
each $j=1,\ldots r$. Since $\Psi$ and $\Phi$ are bounded the
integral converges.  On the other hand, suppose  the integral
converges. Let $\epsilon >0$. Then the integral converges when
integrated over the complement of the set where
$t_i-t_j<\epsilon$, for all $i,j$ such that $1<i<j<r$. Since
$\Phi$ and $\Psi$ are bounded away from zero on such a set it must
be that $-2\nu+(r+1-2j)a<0$ for each $j=1,\ldots r$. This implies
$\nu>(r-1)\frac a 2$.
\end{proof}

We remark that the constant in the last lemma is exactly the same
as the endpoint of the continuous set of parameters, $(\frac a 2
(r-1),\infty)$, for the unitary highest weight modules. Assume
that $\nu>(r-1)\frac{a}{2}$. Let $R_\nu$ be the restriction map
$R_{\nu}:\mathcal{H}_{\nu}(\mathcal{D})\rightarrow C^{\infty
}(\mathcal{D}_{\mathbb{R}})$ given by
\begin{equation}
R_{\nu}f(x)=f(x)h(x)^{\frac{\nu}{2}}.
\end{equation}
As $\nu$ will be fixed most of the time, we will often write $R$ for $R_{\nu}%
$. Consider the restriction of the group action $\pi_{\nu}$ of
$\tilde{G}$ to $\tilde{H}$, where $\tilde{H}$ is the subgroup
corresponding to the Lie algebra $\mathfrak{h}$. Using the method
from \cite{o00}, Lemma 3.4, we can now prove the following lemma.

\begin{lemma}
\label{R} Assume that $\nu>(r-1)\frac{a}{2}$. Then the map $R_\nu$
is a closed densely defined $\tilde{H}$-intertwining operator from
$\mathcal{H}_{\nu }(\mathcal{D})$ into
$L^{2}(\mathcal{D}_{\mathbb{R}},d\eta)$.
\end{lemma}

\begin{proof}
That $R$ intertwines the action of $\tilde{H}$ on
$\mathcal{H}_{\nu }(\mathcal{D})$ and the regular action of
$\tilde{H}$ on $L^{2}(\mathcal{D},d\eta)$ follows by the
transformation properties of $h(x)$. If $F$ is a polynomial,
$RF\in L^{2}(\mathcal{D}_{\mathbb{R}},d\eta)$ by Lemma
\ref{le-sqrint}. Hence $R$ is densely defined.  Let $g\in
C_{c}^{\infty }(\mathcal{D}_{\mathbb{R}})$ and $\epsilon>0$. Then
there exists a polynomial $F$ such that
$||F-h^{-\nu/2}g||_{\infty}<\epsilon/\Vert h^{\nu/2}\Vert$. Hence
\begin{align*}
\Vert RF-g\Vert^{2}  &  =\int_{\mathcal{D}_{\mathbb{R}}}|h(x)^{\nu
}|F(x)-h(x)^{-\nu/2}g(x)|^{2}d\eta(x)\\
&  \leq\frac{\epsilon^{2}}{\Vert h^{\nu/2}\Vert^{2}}\int h(x)^{\nu}%
\,d\eta(x)=\epsilon^{2}%
\end{align*}
Hence $\mathrm{Im}(R)$ is dense in
$L^{2}(\mathcal{D}_{\mathbb{R}},d\eta)$. Finally, $R$ is closed
because point evaluations in $\mathcal{H}_{\nu }(\mathcal{D})$ are
continuous.
\end{proof}

We can now consider the adjoint
$R^{*}:L^{2}(\mathcal{D}_{\mathbb{R}},d\eta)\rightarrow
H_\nu(\mathcal{D})$ as a densely defined operator.

\begin{theorem}
\label{th-main} Assume that $\nu>(r-1)\frac{a}{2}$. Let $f$ be in the domain
of definition of $R_{\nu}^{\ast}$. Then the following holds:

\begin{enumerate}
\item ${\displaystyle
R_\nu R_\nu^{\ast}f(y)=\int_{\mathcal{D}_{\mathbb{R}}}\,\frac{h(y)^{\nu/2}h(x)^{\nu/2}%
}{h(y,x)^{\nu}}\,f(x)\,d\eta(x)\,.}$

\item  If  $g\in H$, then this is the same as
\[
R_{\nu}R_{\nu}^{\ast}f(g\cdot0)=\int_{H}J(h^{-1}g,0)^{\nu/p}f(h)\,d\dot
{h}=D_{\nu}\ast f(g),
\]
where $D_{\nu}(h)=J(h,0)^{\nu/p}$ is in $L^{2}(H/L,d\dot{h})$
and $D_{\nu}\ast f$ stands for the group convolution.

\item  Assume that $\nu>(r-1)a$. Then $D_{\nu}\in L^{1}(\mathcal{D}%
_{\mathbb{R}},d\eta)$.

\item  If $\nu>(r-1)a$. Then $R_{\nu}R_{\nu}^{\ast}:L^{2}(\mathcal{D}%
_{\mathbb{R}},d\eta)\rightarrow
L^{2}(\mathcal{D}_{\mathbb{R}},d\eta)$ is
continuous with norm $\Vert R_{\nu}R_{\nu}^{\ast}\Vert\leq\Vert D_{\nu}%
\Vert_{L^{1}}$.

\item  If $\nu>(r-1)a$. Then $R_{\nu}R_{\nu}^{\ast}:L^{\infty}(\mathcal{D}%
_{\mathbb{R}},d\eta)\rightarrow
L^{\infty}(\mathcal{D}_{\mathbb{R}},d\eta)$ is continuous with
norm $\Vert R_{\nu}R_{\nu}^{\ast}\Vert\leq\Vert D_{\nu
}\Vert_{L^{1}}$.
\end{enumerate}
\end{theorem}

\begin{proof}
(1) Let $f$ be in the domain of definition of $R^{\ast}$. Then $R^{\ast}%
f\in\mathcal{H}_{\nu}(\mathcal{D})$ and for $w\in\mathcal{D}$ we get:
\begin{align*}
R^{\ast}f(w)  &  =(R^{\ast}f,h(\cdot,w)^{-\nu})_{\mathcal{H}_{\nu}%
(\mathcal{D})}\\
&  =(f,R(h(\cdot,w)^{-\nu}))_{L^{2}(\mathcal{D}_{\mathbb{R}})}\\
&
=\int_{\mathcal{D}_{\mathbb{R}}}f(x)h(x)^{\nu/2}h(w,x)^{-\nu}\,d\eta(x),
\end{align*}
where we have used that $\overline{h(z,w)}=h(w,z)$ (c.f. Lemma \ref{le-proph}%
). Thus for $y\in\mathcal{D}_{\mathbb{R}}$ we get:
\[
RR^{\ast}f(y)=\int_{\mathcal{D}_{\mathbb{R}}}f(x)h(y)^{\nu/2}h(x)^{\nu
/2}h(y,x)^{-\nu}\,d\eta(x)
\]
which proves the first statement.

(2) Let $g_1,g_2\in H$. According to Lemma \ref{le-proph} we have
\[
h(g\cdot0)=h(g\cdot0,g\cdot0)=J(g\cdot0,0)^{2/p}.%
\]
Thus
\[
h(g_1\cdot0,g_2\cdot0)=h(g_2 g_2^{-1}g_1\cdot0,g_2\cdot0)=J(g_2,g_2^{-1}g_1\cdot0)^{1/p}%
J(g_2,0)^{1/p},
\]
where we have used that $J(g_2,0)$ and $J(g_1,0)$ are real, and
that $h(z,0)=1$ for all $z$. The cocycle relation
$J(ab,z)=J(a,bz)J(b,z)$ gives
\[
J(g_1,0)=J(g_2(g_2^{-1}g_1),0)=J(g_2,g_2^{-1}g_1\cdot0)J(g_2^{-1}g_1,0)\,.
\]
Hence
\[
J(g_2,g_2^{-1}g_1\cdot0)=J(g_1,0)J(g_2^{-1}g_1,0)^{-1}\,.
\]
Thus the integral kernel for $RR^{\ast}$ becomes
\begin{align*}
\frac{h(g_1\cdot0)^{\nu/2}h(g_2\cdot0)^{\nu/2}}{h(g_1\cdot0,g_2\cdot0)^{\nu}}
&
=\frac{J(g_1,0)^{\nu/p}J(g_2,0)^{\nu/p}}{J(g_1,0)^{\nu/p}J(g_2^{-1}g_1,0)^{-\nu
/p}J(g_2,0)^{\nu/p}}\\
&  =J(g_2^{-1}g_1,0)^{\nu/p}=D_{\nu}(g_2^{-1}g_1)
\end{align*}
As $D_{\nu}(g)=J(g,0)^{\nu/p}=h(g\cdot0)^{\nu/2}$ it follows by
Lemma \ref{le-sqrint} that $D_{\nu}\in L^{2}(H/L,d\dot{h})$.

(3) Notice that $D_{\nu}^{2}=D_{2\nu}$. Hence the claim follows from Lemma
\ref{le-sqrint}.

(4) and (5) are now obvious.
\end{proof}

If $\nu>a(r-1)$ then%
\[
RR^{\ast}1(g\cdot0)=\int_{H}D_{\nu}(g^{-1}h)\,dh=\left| \left|
D_{\nu}\right|  \right|  _{L^{1}}<\infty\,.
\]
Since $ \left|  \left| D_{\nu}\right|  \right|  _{L^{1}}>0$ we
define $c_\nu$ by the relation  $\frac 1 {c_{\nu}}=\left|  \left|
D_{\nu}\right| \right| _{L^{1}}$. The operator $RR^{\ast}$ is
called \textit{the Berezin transform} on
$\mathcal{D}_{\mathbb{R}}$ and $B_{\nu}=c_{\nu}R R^{\ast}$ is
called \textit{the normalized Berezin transform} because
$B_\nu(1)=1$. By theorem (\ref{th-main})
\begin{align}
B_{\nu}f(x)  &  =c_{\nu}\int_{\mathcal{D}_{\mathbb{R}}}\frac{h(x)^{\nu
/2}h(y)^{\nu/2}}{h(x,y)^{\nu}}\,f(y)\,d\eta(y)\nonumber\\
&  =c_{\nu}\int_{\mathcal{D}_{\mathbb{R}}}\frac{h(x)^{\nu/2}h(y)^{\nu/2}%
}{|h(x,y)^{\nu/2}|^{2}}\,f(y)\,d\eta(y)
\end{align}

The space $L^{2}(\mathcal{D}_{\mathbb{R}},d\eta)^{L}$ is
decomposed into a direct integral of principle series
representations of $H$ via the spherical Fourier transform
\begin{equation}
\tilde{f}(\lambda)=\mathcal{F}(f)(\lambda)=\int_{\mathcal{D}_{\mathbb{R}}%
}f(x)\phi_{\lambda}(x)\,d\eta(x)
\end{equation}
or, in terms of the spectral decomposition of commuting
self-adjoint operators, as a direct integral of eigenspaces of the
invariant differential operators. The $L$-invariant eigenfunctions
of $B_\nu$ are precisely the spherical functions
$\phi_{\mathbf{\lambda}}$, $ \lambda \in a^{\ast}_{\mathbb{C}}$.
Thus
\[
B_{\nu}\phi_{\mathbf{\lambda}}=b_{\nu}(\mathbf{\lambda})\phi_{\mathbf{\lambda
}}.
\]
The symbol $b_{\nu}$ is explicitly calculated in \cite{gkz-bere-rbsd}.

\subsection{The generalized Segal-Bargmann transform}

We assume in this section that $\nu>(r-1)\frac{a}{2}$. The
operator $RR^{*}$ is well defined and by definition positive. We
can therefore define $\sqrt{RR^{*}}$. Then there exists a partial
isometry $U_{\nu}$ such that $R^{*}=U_{\nu}\sqrt{RR^{*}}$. To
simplify notation we will often write $U$ for $U_{\nu}$. As
$R=\sqrt{RR^{*}}U^{*}$ and $\mathrm{Im}(R)$ is dense it follows
that $U$ is actually a unitary isomorphism.

\begin{definition}
\label{def-gSBtr} Let $\nu>(r-1)\frac{a}{2}$. The unitary
isomorphism $U_{\nu
}:L^{2}(\mathcal{D}_{\mathbb{R}},d\eta)\rightarrow\mathcal{H}_{\nu
}(\mathcal{D})$ is called the generalized Segal-Bargmann transform
(\cite{OOe96}).
\end{definition}

Let $W=N_{L}(\mathfrak{a})/Z_{L}(\mathfrak{a})$ be the Weyl group
of $\mathfrak{a}$ corresponding to the root system
$\Delta(\mathfrak{h},\mathfrak{a})$ and let
$f\rightarrow\widetilde{f}$ denote the spherical Fourier
transform, Then $f\mapsto\tilde{f}=\mathcal{F}(f)$ extends to an
unitary isomorphism
\begin{equation}
\mathcal{F}:L^{2}(\mathcal{D}_{\mathbb{R}},d\eta)^{L}\rightarrow
L^{2}\left(
\mathfrak{a}^{\ast}/W,\frac{d\lambda}{|c(\lambda)|^{2}}\right)
\simeq L^{2}\left(
\mathfrak{a}^{\ast},\frac{d\lambda}{w|c(\lambda)|^{2}}\right)
^{W}\,,
\end{equation}
where
\[
c(\lambda)=c_{0}\prod_{j<k}\frac{\Gamma(i(\lambda_{j}-\lambda_{k}))}%
{\Gamma((\rho_{j}-\rho_{k}))}\frac{\Gamma(\frac{a}{2}+i(\lambda_{j}%
-\lambda_{k}))}{\Gamma(\frac{a}{2}+\rho_{j}-\rho_{k})}%
\]
is the Harish-Chandra $c$-function, $c_{0}$ is a constant whose
value can be evaluated by using known integral formulas
(see \cite{gkz-jfa} for the case of type C and D domains), and $w$ is the order of the Weyl group $W$.
Combining all of this together we now get the following
proposition.

\begin{proposition} \label{prop-founitary}
Suppose $\nu>(r-1)\frac{a}{2}$. Then
\[
\mathcal{F}\circ
U_\nu^{\ast}:\mathcal{H}_{\nu}(\mathcal{D})^L\rightarrow
L^{2}(\mathfrak{a}^{\ast},\frac{1}{w}|c(\lambda)|^{-2}d\lambda)^{W}%
\]
is a unitary isomorphism.
\end{proposition}

\begin{corollary} \label{cor-orthog}
Let the notation be as above. Then the following hold:

\begin{enumerate}
\item  If $\nu>(r-1)\frac{a}{2}$ then the functions
\[
\mathcal{F}U_{\nu}^{\ast}(\psi_{{\mathbf{m}}})\,,\quad\mathbf{m}%
\in\mathbf{\Lambda}%
\]
form an orthogonal basis for the Hilbert space $L^{2}(\mathfrak a^{\ast
}/W,|c(\lambda)|^{-2}d\lambda)$.

\item  If $\nu>1+a(r-1)$ then the functions
\[
\sqrt{\frac{d_{\mathbf{m}}(\nu)_{\mathbf{m}}}{\left(  \frac{d}{r}\right)
_{\mathbf{m}}}}\,\mathcal{F}U_{\nu}^{\ast}(\psi_{{\mathbf{m}}})\,,\quad
\mathbf{m}\in\mathbf{\Lambda}%
\]
form an orthonormal basis for the Hilbert space $L^{2}(\mathfrak a^{\ast
}/W,|c(\lambda)|^{-2}d\lambda)$.
\end{enumerate}
\end{corollary}

\begin{proof}
This follows Lemma \ref{le-HnuD}, part 6, and Corollary \ref{cor-orthonormal}.
\end{proof}

Our first main goal of this paper is to identify the functions
$\mathcal{F} {U^{*}\psi_{{\mathbf{m}}}}$ and study their
analytical properties.

\section{Generating functions and orthogonality relations for the Branching coefficients}

In this section we derive the orthogonality relations for the
branching coefficients. These results follow somewhat easily from
a general consideration (\cite{gkz-jfa}). In the case where one
considers the branching rules for the tensor product
$\mathcal{H}_{\nu}(\mathcal{D})\otimes\overline{\mathcal{H}_{\nu
}(\mathcal{D})}$ of $\tilde{G}$, considered as the restriction of
the representation of $\tilde{G}\times\tilde{G}$ on the diagonal,
the branching coefficients, also called the Clebsch-Gordan
coefficients, are studied in \cite{gkz-invtoe},
\cite{oz-weyl} and \cite{gkz-tenmin}.
The results here parallel those obtained there. Thus, we will be
brief.

 Define $p_{\mathbf{m}}\in \mathcal{P}(\mathfrak{a}^{\ast})$ by the Rodrigue's type
formula
\begin{equation}
p_{\nu, \mathbf{m}}(\lambda)=p_{\mathbf{m}}(\lambda) :=\|\psi_{\mathbf{m}%
}\|_{\mathcal{F}(V)}^{-2}\psi_{\mathbf{m}}(\partial_{x})(\Delta(e-x^{2}%
)^{-\nu/2}\phi_{\lambda}(x))|_{x=0}.
\end{equation}
We can consider the polynomials $p_{\mathbf{m}}(\lambda)$ as a
generalization of the Hermite polynomials (\cite{gkz-jfa}).

\begin{lemma} \label{le-expansion}
\label{exp} Consider the expansion of $h(x)^{-\nu/2}\phi_{\lambda}%
(x)=\Delta(e-x^{2})^{-\nu/2}\phi_{\lambda}(x)$ in terms of the
$L$-invariant polynomials $\psi_{\mathbf m}$. Then
\begin{equation}
\Delta(e-x^{2})^{-\nu/2}\phi_{\lambda}(x)=\sum_{\mathbf{m}\in\Lambda
}p_{\nu,\mathbf{m}}(\lambda)\psi_{\mathbf{m}}(x) \label{eq-exppm}%
\end{equation}
\end{lemma}

\begin{proof}
Define $P_{\mathbf{m}}(\lambda)$ by
\[
\Delta(e-x^{2})^{-\nu/2}\phi_{\lambda}(x)=\sum_{\mathbf{m}\in\Lambda
}P_{\mathbf{m}}(\lambda)\psi_{\mathbf{m}}(x)\,.
\]
Differentiating both side with respect to
$\psi_{\mathbf{m}}(\partial_{x})$ and setting $x=0$ gives
\[
\Vert\psi_{\mathbf{m}}\Vert_{\mathcal{F}(V)}^{2}p_{\mathbf{m}}(\lambda
)=P_{\mathbf{m}}(\lambda)\left(  \psi_{\mathbf{m}}(\partial_{x})\psi
_{\mathbf{m}}(x)\right)  |_{x=0}=P_{\mathbf{m}}(\lambda)\Vert\psi_{\mathbf{m}%
}\Vert_{\mathcal{F}(V)}^{2}\,.
\]
Hence the claim.
\end{proof}

\begin{lemma}
\label{le-FU} Assume that $\nu>a(r-1)$. Let $F\in\mathcal{H}_{\nu}%
(\mathcal{D})^L$. Then
\[
\mathcal{F}(\sqrt{R_{\nu}R_\nu^{\ast}}U_{\nu}^{\ast}(F))(\lambda)=c_{\nu}^{-1/2}b_{\nu
}(\lambda)^{1/2}\,\mathcal{F}(U_{\nu}^{\ast}(F))(\lambda)\,.
\]
\end{lemma}

\begin{proof}
Recall that $RR^{\ast}=D_{\nu}\ast$. Hence, for all $f\in
L^{2}(\mathcal{D}_{\mathbb{R}},d\eta)^L$,
\[
\mathcal{F}(R R^{\ast}f)(\lambda)=\mathcal{F}(D_{\nu})(\lambda
)\mathcal{F}(f)(\lambda)\,.
\]
Furthermore,
\[
\mathcal{F}(D_{\nu})(\lambda)=\int
D_{\nu}(h)\phi_{\lambda}(h)\,d\dot
{h}=c_{\nu}^{-1}B_{\nu}(\phi_{\lambda})(1)=c_{\nu}^{-1}b_{\nu}(\lambda),
\]
since $\phi_{\lambda}(1)=1$. Consequently,
\[
\mathcal{F}(\sqrt{ RR^{\ast}}f)(\lambda)=c_{\nu}^{-\frac 12}\sqrt{b_{\nu}(\lambda)}%
\,\tilde{f}(\lambda)\,.
\]
Applying this to a function $f=U^{\ast}F$, $F\in\mathcal{H}_{\nu}%
(\mathcal{D})$, gives  the lemma.
\end{proof}

The next proposition states that the polynomials $p_{\mathbf{m}}(\lambda)$ can
be also obtained via the Segal-Bargman transform of $\psi_{\mathbf{m}}$.

\begin{proposition} \label{prop-poly}
Assume that $\nu>a(r-1)$. Then
\[
\mathcal{F}(U_{\nu}^{\ast}\psi_{\mathbf{m}})(\lambda)=c_{\nu}^{-\frac 12}%
\sqrt{b_{\nu}(\lambda)}\Vert\psi_{\mathbf{m}}\Vert_{\nu}^{2}\,\,p_{\nu,\mathbf{m}%
}(\lambda)
\]
\end{proposition}

\begin{proof}
Since  $\nu>a(r-1)$, we have
\[
\int_{D}h(x)^{\frac{\nu}{2}}d\eta(x)<\infty.
\]
Let $C(\rho)$ be the convex hull of $W\cdot\rho$. By the
Helgason-Johnson theorem ( Theorem 8.1 of \cite{Helg}, p. 458, see
also \cite{HJ} ), if $\lambda\in\mathfrak a^{\ast}+iC(\rho)$ then
$\phi_{\lambda}$ is a bounded function on
$\mathcal{D}_{\mathbb{R}}$. Thus the integral defining
$B_{\nu}\phi _{\mathbf{\lambda}}$ is absolutely convergent and
\[
B_{\nu}\phi_{\lambda}(x)=c_{\nu}\int_{\mathcal{D}_{\mathbb{R}}}\frac
{h(x)^{\nu/2}h(y)^{\nu/2}}{h(x,y)^{\nu}}\phi_{\lambda}(y)\,d\eta(y)=b_{\nu
}(\lambda)\phi_{\lambda}(x)\,.
\]
We divide by $c_{\nu}h(x)^{\nu/2}$. Furthermore, extend $\Vert
\psi_{\mathbf{m}}\Vert_\nu^{-1}\psi_{\mathbf{m}}$, $\mathbf{m}\in\mathbf{\Lambda}%
$, to an orthonormal basis $F_{\mathbf{n}}$ of $\mathcal{H}_{\nu}%
(\mathcal{D})$. Then
\[
h(z,w)^{-\nu}=\sum_{\mathbf{n}}F_{\mathbf{n}}(z)\overline{F_{\mathbf{n}}%
(w)}\,.
\]
This gives:
\[
c_{\nu}^{-1}b_{\nu}(\lambda)h(x)^{-\nu/2}\phi_{\lambda}(x)=\int_{\mathcal{D}%
_{\mathbb{R}}}\sum_{\mathbf{n}}h(y)^{\nu/2}F_{\mathbf{n}}(y)\phi_{\lambda
}(y)\overline{F_{\mathbf{n}}(x)}\,d\eta(y)\,.
\]
We now let $\psi_{\mathbf{m}}(\partial_{x})\vert_{x=0}$ act on
both sides. For the left-hand side we use Lemma
\ref{le-expansion}. For the right-hand side we use the fact that
$\Vert
\psi_{\mathbf{m}}\Vert_\nu^{-1}\psi_{\mathbf{m}}(\partial_{x})\overline{F_{\mathbf{n}}(x)}%
|_{x=0}=\delta_{\mathbf{m}\mathbf{n}}\frac{\Vert\psi_{\mathbf{m}}\Vert_{\mathcal{F}(V)}^{2}}{\Vert\psi_{\mathbf{m}%
}\Vert_{\nu}^{2}}$, when $F_{\mathbf m}$ is one of the $\Vert
\psi_{\mathbf{m}}\Vert_\nu^{-1}\psi_{\mathbf{m}}$, and $0$
otherwise. We thus get
\[\Vert\psi_{\mathbf{m}}\Vert_{\mathcal{F}%
(V)}^{2}c_{\nu}^{-1}b_{\nu}(\lambda)p_{\mathbf{m}}(\lambda)=
\frac{\Vert\psi_{\mathbf{m}}\Vert_{\mathcal{F}(V)}^{2}}{\Vert\psi_{\mathbf{m}%
}\Vert_{\nu}^{2}}\int_{\mathcal{D}_{\mathbf{R}}}h(y)^{\nu/2}\psi_{\mathbf{m}%
}(y)\phi_{\lambda}(y)\,d\eta(y).
\]
That is,
\[
\mathcal{F}(R\psi_{\mathbf{m}})(\lambda)=\Vert\psi_{\mathbf{m}}\Vert_{\nu}%
^{2}c_{\nu}^{-1}b_{\nu}(\lambda)p_{\mathbf{m}}(\lambda)
\]
The claim follows now from Lemma \ref{le-FU} as
$R=\sqrt{RR^{\ast}}\,U^{\ast}$.
\end{proof}

\section{Recurrence formulas for the Branching coefficients}
Propositions (\ref{prop-founitary}) and (\ref{prop-poly}) together
imply that when the measure $b(\lambda)\vert c(\lambda)\vert^{-2}
d\lambda$ is associated with  $\mathfrak{a}^{\ast}/W$ then the
polynomials $p_{\mathbf{n}}(\lambda), \; \mathbf{n}\in \Lambda$,
are orthogonal. When reference is made to the orthogonality of
these polynomials it is understood to be with respect to this
measure.
 In this section and the next we will derive the recurrence and
difference formulas for the orthogonal polynomials
$p_{\mathbf{n}}(\lambda)$. In most cases, our method will be to
prove our formula first for $\lambda$ in a certain integral cone
where we can easily deal with the equations. To pass to the
general formula we then need the following elementary result,
which can be easily proved by induction. Recall the lattice
$\mathbf{\Lambda}$ defined in (\ref{eq-lambda}).

\begin{lemma}
\label{reduc} Suppose $p_{1}(\lambda)$ and $p_{2}(\lambda)$ are two
polynomials in $\mathfrak{a}_{\mathbb{C}}^{\ast}$. Let $\mu\in\mathfrak{a}%
_{\mathbb{C}}^{\ast}$. Suppose that $p_{1}(\lambda)=p_{2}(\lambda)$ for all
$\lambda=\mathbf{n}+\mu$, $\mathbf{n}\in\mathbf{\Lambda}$. Then $p_{1}%
(\lambda)=p_{2}(\lambda)$ for all
$\lambda\in\mathfrak{a}_{\mathbb{C}}^{\ast}$.
\end{lemma}

To simplify certain arguments we assume, in this and next section, that $\nu$
is an even positive integer. To state our result we need the binomial
coefficient
\[
\binom{\mathbf{n}}{\mathbf{n}-\gamma_{k}}=(n_{k}+\frac{a}{2}(r-k))\prod_{j\neq
k}\frac{n_{k}-n_{j}+\frac{a}{2}(j-k-1)}{n_{k}-n_{j}+\frac{a}{2}(j-k)}\,,
\]
(c.f. \cite{Lassalle-binom}). Note that $\binom{\mathbf{n}}{\mathbf{n}%
-\gamma_{k}}$ is actually a rational function of $\mathbf{n}$ and
can be defined on all $\mathbb{C}^{n}$ with singularities on some
lower dimensional hyperplanes. We let
\[
c_{\mathbf{n}}(k)=\prod_{j\neq k}\frac{n_{j}-n_{k}-\frac{a}{2}(j+1-k)}%
{n_{j}-n_{k}-\frac{a}{2}(j-k)}.%
\]
The following theorem  is essentially proved in \cite{oz-duke},
although the holomorphic discrete series $\pi_{\nu}$ is realized
there in a degenerate principal series representation. We will not
reproduce the same proof here.
 In the case of $G=\mathrm{SU}%
(1,1)$ we have $\psi_{m}(z)=z^{m}$. This theorem and the lemma
that follows are
generalizations of the simple facts that%
\begin{align*}
\pi_{\nu}\left(  -\left(
\begin{array}
[c]{cc}%
0 & 1\\
1 & 0
\end{array}
\right)  \right)  z^{m}  &  =-(\nu+m)z^{m-1}+mz^{m-1}\\
\pi_{\nu}\left(  -\left(
\begin{array}
[c]{cc}%
1 & 0\\
0 & -1
\end{array}
\right)  \right)  z^{m}  &  =(\nu+2m)z^{m},%
\end{align*}
which can be verified by a simple calculation.

\begin{theorem} Recall that $\xi=\mathrm{Ad}(\gamma)(-2Z_{0})$.
\label{xi-act} In the bounded realization $\mathcal{H}_{\nu}(\mathcal{D})$ of
the representation of $\mathfrak g_{\mathbb{C}}$,
\begin{equation}
\pi_{\nu}(-\xi)\psi_{\mathbf{n}}=\sum_{j=1}^{r}\binom{\mathbf{n}}%
{\mathbf{n}-\gamma_{j}}\psi_{{\mathbf{n}}-\gamma_{j}}-\sum_{j=1}^{r}%
(\nu+{n}_{j}-\frac{a}{2}(j-1))c_{\mathbf{n}}(j)\psi_{{\mathbf{n}%
}+\gamma_{j}} \label{eq-xiact}%
\end{equation}
\end{theorem}

\begin{remark}
As is proved in \cite{oz-duke} one may use the above formula to
derive the results of Faraut-Koranyi \cite{FK} on the norm of
$\phi_{\mathbf{n}}$ in the Hilbert space
$\mathcal{H}_{\nu}(\mathcal{D})$. Conversely, knowing that
$\pi_{\nu}(\xi)\phi_{\mathbf{\lambda}}$ is of the above form, it
is possible also to find recursively the coefficients by using the
results in \cite{FK}, the dimension formula in \cite{UP-tams} and
the fact that $\pi_{\nu}(\xi)$ is a skew-symmetric operator.
\end{remark}

The element $2Z_{0}=D(e,\bar{e})$ is in the center of $\mathfrak
k$, and thus acts on each $\psi_{\mathbf{n}}$ by a scalar. A
routine calculation gives the following lemma:

\begin{lemma}
\label{le-actz0} We have
\begin{equation}
\pi_{\nu}(-2Z_{0})\psi_{\mathbf{m}}=(r\nu+2|\mathbf{m}|)\psi_{\mathbf{m}}.
\label{eq-actz0}%
\end{equation}
\end{lemma}

We can draw several important relations from these facts. Recall
first from
Theorem \ref{th-2.5}, part 5 that%
\[
q_{\mathbf{m},\nu}(z)=\Delta(z+e)^{-\nu}\psi_{\mathbf{m}}\left(
(z-e)(z+e)^{-1}\right)
=2^{r\nu/2}\,\pi_{\nu}(\gamma)\psi_{\mathbf{m}}(z)\,.
\]
Recall that $\pi_{\nu}(\gamma)$ is well defined and that
\[
\pi_{\nu}(\gamma)\pi_{\nu}(X)=\pi_{\nu}(\mathrm{Ad}(\gamma)X)\pi_{\nu}%
(\gamma)
\]
for all $X\in\mathfrak{g}_{\mathbb{C}}$. As
$\mathrm{Ad}(\gamma^{-1})(2Z_{0})=\xi$ and
$\mathrm{Ad}(\gamma^{-1})(\xi)=-2Z_{0}$, by Lemma \ref{le-xi}, we
get by applying the Cayley transform to (\ref{eq-xiact}) and
(\ref{eq-actz0}):

\begin{lemma}
\label{le-rerelqmnu}Let $q_{\mathbf{m},\nu}(z)=\Delta(z+e)^{-\nu}%
\psi_{\mathbf{m}}\left(  (z-e)(z+e)^{-1}\right)  $ then the following holds:
\begin{enumerate}
\item $\pi_{\nu}(\xi)q_{\mathbf{m},\nu}=(r\nu+2\left|\mathbf{m}\right|
)q_{\mathbf{m},\nu}$.
\item $\pi_{\nu}(-2Z_{0})q_{\mathbf{m},\nu}=\sum_{j=1}^{r}\binom{\mathbf{m}%
}{\mathbf{m}-\gamma_{j}}q_{{\mathbf{m}}-\gamma_{j},\nu}-\sum_{j=1}^{r}%
(\nu+{n}_{j}-\frac{a}{2}(j-1)))c_{\mathbf{m}}(j)q_{{\mathbf{m}}%
+\gamma_{j},\nu}$.
\end{enumerate}
\end{lemma}

The following theorem gives the recursion relations for the
polynomials $p_{\nu,\mathbf{m}}$.

\begin{theorem}\label{thm-recrel}
The following recurrence formula holds:
\[
(2\sum_{j=1}^{r}(i\lambda_{j}+\rho_{j}))p_{\nu,\mathbf{m}}(\lambda)=\sum_{j=1}%
^{r}\binom{\mathbf{m}+\gamma_{j}}{\mathbf{m}}p_{\nu,\mathbf{m}+\gamma_{j}}%
(\lambda)-(\nu+m_{j}-1-\frac{a}{2}(j-1))c_{\mathbf{m}-\gamma_{j}%
}(j)p_{\nu,\mathbf{m}-\gamma_{j}}(\lambda).
\]
\end{theorem}

\begin{proof}
We prove the result for
\begin{equation}
\lambda=-i(\mathbf{m}+\frac{\nu}{2}-\rho) \label{lam-m-rel},%
\end{equation}
where $\mathbf{m}\in\Lambda$ and $\frac{\nu}{2}$ is viewed as $\frac{\nu}%
{2}\sum_{j=1}^{r}\beta_{j}$. The general result for $\lambda\in\mathbb{C}%
^{r}=\mathfrak a_{\mathbb{C}}^{\ast}$ follows from Lemma
\ref{reduc}. As
$\mathrm{Ad}(\gamma)\xi=-\mathrm{Ad}(\gamma)^{-1}\xi=2Z_{0}$ we
get from Lemma \ref{le-rerelqmnu}:
\[
\pi_{\nu}(\xi)\pi_{\nu}(\gamma^{-1})\psi_{\mathbf{m}}(x)=-(r\nu+2|\mathbf{m}|)
\pi_{\nu}(\gamma^{-1})\psi_{\mathbf{m}}(x).
\]
Now simplify this expression before actually performing the differentiation:
\begin{eqnarray}
\pi_{\nu}(\gamma^{-1})\psi_{\mathbf{m}}(x)&=&2^{r\nu/2}\psi_{\mathbf{m}%
}\left(  \frac{e+x}{e-x}\right)  \Delta(e-x)^{-\nu}\\
&  =&2^{r\nu/2}\psi_{{\mathbf{m}}+\frac{\nu}{2}}\left(  \frac{e+x}%
{e-x}\right)  \Delta\left(  \frac{e+x}{e-x}\right)  ^{-\frac{\nu}{2}}%
\Delta(e-x)^{-\nu}\nonumber\\
&  =&2^{r\nu/2}\psi_{\mathbf{m}+\frac{\nu}{2}}\left(
\frac{e+x}{e-x}\right)
\Delta(e-x^{2})^{-\frac{\nu}{2}}\nonumber\\
&=&2^{r\nu/2}\phi_{\lambda}(x)\Delta(e-x^{2})^{-\frac{\nu}{2}}\nonumber\,.
\label{eq-gammapsi}
\end{eqnarray}
Here we have used that
\begin{align*}
\Delta\left(  \frac{e+x}{e-x}\right)  ^{-\nu/2}\Delta(e-x)^{-\nu}  &
=\Delta(e+x)^{-\nu/2}\Delta(e-x)^{\nu/2}\Delta(e-x)^{-\nu}\\
&  =\Delta(e+x)^{-\nu/2}\Delta(e-x)^{-\nu/2}\\
&  =\Delta(e-x^{2})^{-\nu/2}\,.
\end{align*}
Hence
\begin{equation}
\pi_{\nu}(\xi)(\Delta(e-x^{2})^{-\nu/2}\phi_{\lambda}(x))=-(r\nu
+2|\mathbf{m}|)\Delta(e-x^{2})^{-\nu/2}\phi_{\lambda}(x)\,.
\label{eq-5.8}%
\end{equation}
On the other hand, by Lemma \ref{exp}:
\begin{equation}
\phi_{\lambda}(x)\Delta(e-x^{2})^{-\frac{\nu}{2}}=\sum_{\mathbf{n}\in\Lambda
}p_{\mathbf{n}}(\lambda)\psi_{\mathbf{n}}(x). \label{exp-shift}%
\end{equation}
Thus (\ref{eq-5.8}) reads
\begin{equation}
\pi_{\nu}(\xi)\sum_{\mathbf{n}\in\Lambda}p_{\mathbf{n}}(\lambda)\psi
_{\mathbf{n}}(x)=-(r\nu+2|\mathbf{m}|)\sum_{\mathbf{n}\in\Lambda}%
p_{\mathbf{n}}(\lambda)\psi_{\mathbf{n}}(x)\,. \label{eq-5.9}%
\end{equation}
Notice that (\ref{exp-shift}) is the power series expansion of an
analytic function and the operator $\pi_{\nu}(\xi)$ is a
differential operator; it commutes with the summation. We thus
have, in a neighborhood of $x=0$,
\begin{equation}
\sum_{\mathbf{n}\in\Lambda}p_{\mathbf{n}}(\lambda)\pi_{\nu}(\xi)\psi
_{\mathbf{n}}(x)=-(r\nu+2|\mathbf{m}|)\sum_{\mathbf{n}\in\Lambda}%
p_{\mathbf{n}}(\lambda)\psi_{\mathbf{n}}(x). \label{eqone}%
\end{equation}
The left hand side, in view of Theorem \ref{xi-act}, is
\begin{align}
&  \sum_{\mathbf{n}\in\Lambda}p_{\mathbf{n}}(\lambda)\left(  \sum_{j=1}%
^{r}\binom{\mathbf{n}}{\mathbf{n}-\gamma_{j}}\psi_{\mathbf{n}-\gamma_{j}%
}(x)-\sum_{j=1}^{r}(\nu+n_{j}-\frac{a}{2}(j-1))c_{\mathbf{n}}(j)\psi
_{\mathbf{n}+\gamma_{j}}(x)\right) \nonumber\\
&  \qquad=\sum_{\mathbf{n}\in\Lambda}\psi_{\mathbf{n}}(x)\left(  \sum
_{j=1}^{r}\binom{\mathbf{n}+\gamma_{j}}{\mathbf{n}}p_{\mathbf{n}+\gamma_{j}%
}(\lambda)-(\nu+n_{j}-1-\frac{a}{2}(j-1))c_{\mathbf{n}-\gamma_{j}%
}(j)p_{\mathbf{n}-\gamma_{j}}(\lambda)\right)  \label{eqtwo}%
\end{align}
Equating the coefficients of $\psi_{\mathbf{n}}(x)$ in (\ref{eqone}) and
(\ref{eqtwo}) we get
\[
(r\nu+2|\mathbf{m}|)p_{\mathbf{n}}(\lambda)=\sum_{j=1}^{r}\binom
{\mathbf{n}+\gamma_{j}}{\mathbf{n}}p_{\mathbf{n}+\gamma_{j}}(\lambda
)-(\nu+n_{j}-1-\frac{a}{2}(j-1))c_{\mathbf{n}-\gamma_{j}}(j)p_{\mathbf{n}%
-\gamma_{j}}(\lambda).
\]
The relation (\ref{lam-m-rel}) implies that $(r\nu+2|\mathbf{m}%
|)=2|i\lambda+\rho|=2\sum_{j=1}^{r}(i\lambda_{j}+\rho_{j})$. This
finishes the proof.
\end{proof}

\begin{example}
If $\mathcal{D}$ is the unit disk, then we can take
$\mathcal{D}_{\mathbb{R}}$ as the unit interval $(-1,1)$ on the
real line. The spherical function on the unit disk is
$\phi_{\lambda}(x)=(\frac{1+x}{1-x})^{i\lambda}$ and the expansion
(\ref{eq-exppm}) reads
\begin{equation}
(1-x^{2})^{-\frac{\nu}{2}}\left(  \frac{1+x}{1-x}\right)
^{i\lambda
}=(1-x)^{-\frac{\nu}{2}-i\lambda}(1+x)^{-\frac{\nu}{2}+i\lambda}=\sum
_{n=0}^{\infty}p_{n,\nu}(\lambda)x^{n} \label{gen-func-disk}%
\end{equation}
with
\[
p_{n,\nu}(\lambda)=(\frac{\nu}{2}+i\lambda)_n
{}_{2}F_{1}(\frac{\nu}{2}%
-i\lambda,-n;-\frac{\nu}{2}-i\lambda-n+1,1)
\]
being the Meixner-Pollacyck polynomials, (c.f. \cite{Basu-Wolf}).
The action of
$\xi=%
\begin{pmatrix}
0 & 1\\
1 & 0
\end{pmatrix}
$ on functions on $(-1,1)$ is
\[
\pi_{\nu}(\xi)f(x)=\nu xf(x)-(1-x^{2})f'(x).
\]
Writing the function in equation (\ref{gen-func-disk}) as
$G_{\nu,\lambda}(x)$ we have
\[
\pi_{\nu}(\xi)G_{\nu,\lambda}=(-2i\lambda)G_{\nu,\lambda}.%
\]
This can be proved easily by a direct computation. It exemplifies
equation (\ref{eq-5.8}) (Here, $\rho=0$ and so $i\lambda= m +
\frac \nu 2$). It then implies the recurrence relation,
\[
(2i\lambda)p_{\nu,n}(\lambda)=(n+1)p_{\nu,n+1}(\lambda)-(\nu+n-1)p_{\nu
,n-1}(\lambda),
\]
and this coincides with Theorem \ref{thm-recrel}.
\end{example}

\begin{remark}
While deriving the recurrence formula we have extended the action
of $\mathfrak g_{\mathbb{C}}$ and $\gamma$ on
$\mathcal{H}_{\nu}(T(\Omega))$ (or
$\mathcal{H}_{\nu}(\mathcal{D})$) to the space of meromorphic
functions on $V$. Indeed the operator $\pi_{\nu}(\gamma^{-1})$ is
up to a constant, c.f., Lemma \ref{le-HnuD}, a unitary operator
from $\mathcal{H}_{\nu}(T(\Omega))$ onto
$\mathcal{H}_{\nu}(\mathcal{D})$ so it is initially defined on
$\mathcal{H}_{\nu}(T(\Omega))$. However, in our formulas the
action of $\pi_{\nu}(\gamma^{-1})$ on $\psi_{\lambda}$ for
$\lambda=\mathbf{m}$ is viewed as the extended action since
$\psi_{\lambda}$ is not an element in
$\mathcal{H}_{\nu}(T(\Omega))$. It suggests some interesting
applications of the idea of extending holomorphic functions on a
domain to meromorphic functions to a larger domain.
\end{remark}

\section{Difference formulas for the Branching coefficients}

In this section we state and prove the difference equation for the
polynomials $p_{\nu,{\mathbf{m}}}(\lambda)$

\begin{theorem}\label{th-diff}
The polynomials $p_{\nu,{\mathbf{m}}}(\lambda)$ satisfy the
following
difference equation%
\begin{align*}
-(r\nu &+ 2|{\mathbf{n}}|)p_{\nu,{\mathbf{n}}}(\lambda) = \\
  &\sum_{j=1}^{r}%
\binom{i\lambda+{\mathbf{\rho}}-\frac{\nu}{2}}{i\lambda+\rho-\frac{\nu}%
{2}-\gamma_{j}}p_{\nu,{\mathbf{n}}}(\lambda+i\gamma_{j})
  -\sum_{j=1}^{r}(\frac{\nu}{2}+i\lambda_{j}+\rho_{j}-\frac{a}%
{2}(j-1)))c_{i\lambda+\rho-\frac \nu 2
}(j)p_{\nu,{\mathbf{n}}}(\lambda-i\gamma_{j}).
\end{align*}
\end{theorem}

\begin{proof}
As in the proof of Theorem \ref{thm-recrel} it suffices to prove
the theorem for those $\lambda$ satisfying  $i\lambda
+\rho={\mathbf{m}}+\nu/2$. Equations (\ref{eq-gammapsi}) and
(\ref{exp-shift}) in the proof above combine to give
\begin{equation} \label{eq-6.1}
2^{\frac{-r\nu}{2}}\pi_\nu(\gamma^{-1})\psi_{\mathbf{m}}=
\sum_{\mathbf{n}\in\Lambda}p_{\mathbf{n}}(\lambda)\psi_{\mathbf{n}}.
\end{equation}
Let the operator $\pi_{\nu}(-2Z_{\circ})$ act on both sides. For
the left-hand side we use Theorem \ref{thm-recrel},
Lemma \ref{le-xi}, and equation (\ref{eq-6.1}) applied to the
case $i(\lambda\pm i\gamma_j)+\rho=(\mathbf{m}\mp\gamma_j) +\frac
\nu 2$ to obtain
\begin{eqnarray*}
LHS&=&2^{\frac
{-r\nu}{2}}\pi_\nu(-2Z_{\circ})\pi_\nu(\gamma^{-1})\psi_{\mathbf{m}}\\
&=&2^{\frac{-r\nu}{2}}\pi_\nu(\gamma^{-1})\pi_\nu(\xi)\psi_{\mathbf{m}}\\
&=&-2^{\frac{-r\nu}{2}}\pi_\nu(\gamma^{-1})
(\sum_{j=1}^{r}\binom{\mathbf{m}} {\mathbf{m}-\gamma_{j}}
\psi_{{\mathbf{m}}-\gamma_{j}}-\sum_{j=1}^{r}
(\nu+{m}_{j}-\frac{a}{2}(j-1))c_{\mathbf{m}}(j)\psi_{{\mathbf{m}}+\gamma_{j}})\\
&=&-\sum_{\mathbf{n}\in\Lambda}\left(\sum_{j=1}^r
\binom{\mathbf{m}}{\mathbf{m}-\gamma_j}p_{\mathbf{n}}(\lambda +
i\gamma_j)-\sum_{j=1}^r(\nu+m_j-\frac a
2(j-1))c_{\mathbf{m}}(j)p_{\mathbf{n}}(\lambda-i\gamma_j)
\right)\psi_{\mathbf{n}}.
\end{eqnarray*}
For the right-hand side we obtain by Lemma \ref{le-actz0}
\begin{eqnarray*}
RHS&=&\pi_\nu(-2Z_{\circ})\sum_{\mathbf{n}\in\Lambda}p_{\mathbf{n}}(\lambda)\psi_{\mathbf{n}}\\
&=&\sum_{\mathbf{n}\in\Lambda}(r\nu +
2\vert\mathbf{n}\vert)p_{\mathbf{n}}(\lambda)\psi_{\mathbf{n}}.
\end{eqnarray*}
The proof is completed by equating the coefficients, rewriting
each occurrence of $\mathbf{m}$ in terms of $\lambda$, and then
applying Lemma \ref{reduc}.
\end{proof}

\begin{example}
We continue the Example 5.7. The above difference relation
can be proved by simple however tricky computations, which in turn reveal
the advantage of using the representation theoretic method.
Using the notation there,
write $G_{\lambda}(x)=G_{\nu, \lambda}(x)$
the generating function. Differentiating the expansion
we get
\begin{equation}
\label{ex6.2-1}
2x\frac{d}{dx}G_{\lambda}(x)=
\sum_{n=0}^\infty 2n p_{n, \nu}(\lam)x^n
\end{equation}
On the other hand, differentiating
the formula for $G_{\lambda}(x)$ results in
\begin{equation}
\label{ex6.2-2}
2x\frac{d}{dx}G_{\lambda}(x)
=(\frac \nu 2 +i\lambda)\frac{2x}{1+x}G_{\lam-i}(x)
+(-\frac \nu 2 +i\lambda)\frac{2x}{1-x}G_{\lam+i}(x).
\end{equation}
We observe  that
$$
G_{\lam}(x)=\frac{1-x}{1+x}G_{\lam-i}(x)=\frac{1+x}{1-x}G_{\lam+i}(x),
$$
which then imply that
\begin{equation}
\label{ex6.2-3}
\nu\, G_{\lam}(x)=
(\frac \nu 2 +i\lam)
G_{\lam}(x) +
(\frac \nu 2 -i\lam)
G_{\lam}(x)
=(\frac \nu 2 +i\lam)
 \frac{1-x}{1+x} G_{\lam-i}(x)
+(\frac \nu 2 -i\lam)
 \frac{1+x}{1-x} G_{\lam-i}(x)
\end{equation}
Summing the equations
(\ref{ex6.2-2})
and (\ref{ex6.2-3}) gives
$$(2x\frac{d}{dx}+\nu)G_{\lambda}(x)
=(\frac \nu 2 +i\lam)G_{\lam-i}(x)
+(\frac \nu 2 -i\lam)G_{\lam+i}(x)
$$
and consequently
$$
(2n+\nu)p_{\nu, n}(\lam)
=(\frac \nu 2 +i\lam)p_{\nu, n}(\lam-i)
+(\frac \nu 2 -i\lam)p_{\nu, n}(\lam+i)
$$
or
$$
-(2n+\nu)p_{\nu, n}(\lam)
=(-\frac \nu 2 +i\lam)p_{\nu, n}(\lam+i)
-(\frac \nu 2 +i\lam)p_{\nu, n}(\lam-i)
$$
which coincides with Theorem 6.1. The  proof
of Theorem 6.1 is thus conceptually clearer.
\end{example}

\section{The restriction Principle in the unbounded realization}

In this section we discuss the application of the restriction
principle to the unbounded realization of $G/K$. In particular, we
use this to introduce the Laguerre functions of the cone $\Omega$
and derive some relations they satisfy. For a function $F$ defined
on the Siegel domain $T(\Omega )=\Omega+iJ$ let
\[
R_{\nu}F(x)=RF(x)=F(x)\,,
\]
where $x\in\Omega$. Since the functions in $\mathcal{H}_{\nu}^{2}(T(\Omega))$
are holomorphic it follows as before that $R$ is injective on $\mathcal{H}%
_{\nu}^{2}(T(\Omega))$. For $y\in\Omega$, let $k_{y}(z)=K_{\nu}(z,y)=K(z,y)$,
where $z\in T(\Omega)$. Then $k_{y}\in\mathcal{H}_{\nu}(T(\Omega))$ and%
\[
Rk_{y}(x)=\Gamma_{\Omega}(\nu)\Delta(x+y)^{-\nu}\,.
\]

\begin{lemma}
The linear span of $\{k_{y}\mid y\in\Omega\}$ is dense in $\mathcal{H}_{\nu
}(T(\Omega))$.
\end{lemma}

\begin{proof}
Assume that $f\in\mathcal{H}_{\nu}(T(\Omega))$ is perpendicular to all $k_{y}%
$, $y\in\Omega$. Then $f(y)=0$, for all $y\in\Omega$. As $f$ is
holomorphic it follows that $f=0$.
\end{proof}

Recall the \textit{Beta-function} for $\Omega$:%
\[
B_{\Omega}(\nu,\mu)=\int_{\Omega}\Delta(x+e)^{-\nu-\mu}\Delta(x)^{\nu
-d/r}\,dx\,,
\]
which is finite for $\nu$ and $\mu$ having real parts greater than
$ (r-1)\frac a 2$ (\cite{FK-book}, p 141).
\begin{lemma}
Let $y\in\Omega$. Then $Rk_{y}\in L_{\nu}^{2}(\Omega)$, for all $\nu
>(r-1)\frac{a}{2}$. In fact%
\[
\left|  \left|  Rk_{y}\right|  \right|  =\frac{\Gamma_{\Omega}(\nu)}%
{\Delta(y)^{\nu/2}}\sqrt{B_{\Omega}(\nu,\nu)}.
\]
\end{lemma}

\begin{proof}
Let $y\in\Omega$. Then, performing the change of variables
$x\mapsto Q(y^{\frac 12})x$ where $Q$ is the Jordan quadratic
operator in the section 2.1,
\begin{align*}
\left\|  Rk_{y}\right\|  ^{2}  &  ={\Gamma_{\Omega}^{2}(\nu)}\int_{\Omega
}\frac{\Delta(x)^{\nu-d/r}}{\Delta(x+y)^{2\nu}}\,dx\\
&  =\Gamma_{\Omega}^{2}(\nu)\frac{1}{\Delta(y)^{\nu}}\int_{\Omega}\frac
{\Delta(x)^{\nu-d/r}}{\Delta(x+e)^{2\nu}}\,dx\quad
\text{ ($x\rightarrow
Q(y^{\frac{1}{2}}) x$)}\\
&  =\frac{\Gamma_{\Omega}^{2}(\nu)}{\Delta(y)^{\nu}}B_{\Omega}(\nu,\nu
)<\infty,
\end{align*}
finishing the proof.
\end{proof}

We need to  distinguish the Laplace transform $\mathcal L_\nu f$
as mapping on functions on $\Ome$ to itself and as mapping
functions on $\Ome$ to holomorphic functions on $T(\Ome)$. We
write the former as $\mathcal L_\nu^\Ome f$, which is the
restriction of the latter.

\begin{lemma} Let $\nu >(r-1)\frac{a}{2}$. Then
the set $\left\{  Rk_{y}\mid y\in\Omega\right\}  $ is dense in $L_{\nu}%
^{2}(\Omega)$.
\end{lemma}

\begin{proof}
Let $f\in L_{\nu}^{2}(\Omega)$ and suppose $f$ is orthogonal to all $Rk_{y}$,
$y\in\Omega$. Then
\begin{align*}
0  &  =\left(  f\mid Rk_{y}\right) \\
&  ={\Gamma_{\Omega}(\nu)}\int_{\Omega}f(x)\Delta(x+y)^{-\nu}\Delta
(x)^{\nu-d/r}dx\\
&  =\int_{\Omega}f(x)\mathcal{L}_{\nu}(e^{-(y|\cdot)})(x)\Delta(x)^{\nu
-d/r}dx\\
&  =\int_{\Omega}e^{-(y|x)}\mathcal{L}^\Ome_{\nu}(f)(x)\Delta(x)^{\nu-d/r}\,dx\\
&  =\mathcal{L}^\Ome_{\nu}(\mathcal{L}^\Ome_{\nu}(f))(y)\,.
\end{align*}
{} From this and the injectivity of the Laplace transform it follows that
$f=0$.
\end{proof}

By the previous lemmas it follows that the restriction map
\[
R:\mathcal{H}_{\nu}(T(\Omega))\rightarrow L_{\nu}^{2}(\Omega)
\]
is injective, densely defined and has dense range.  Since
convergence in $\mathcal{H}_{\nu}(T(\Omega))$ implies uniform
convergence on compact sets it follows easily that $R$  is closed.
We can thus polarize $R^{\ast}$, $R^{\ast}=U|R^{\ast}|$, and
obtain a unitary map $U$
 from $L_{\nu}^{2}(\Omega)$ onto
$\mathcal{H}_{\nu}(T(\Omega))$. We now prove that
$U=\mathcal{L}_\nu.$

\begin{theorem}\label{th-Laplace}
\label{Th-7.4} Suppose $\nu >\frac a2(r-1)$. The polar
decomposition of $R^\ast$ is given by $R^\ast=\mathcal L_\nu
\mathcal L_\nu^\Ome$, where  the operator $\mathcal L_\nu$ extends
to a unitary operator from $L_{\nu}^{2}(\Omega)$ onto
$\mathcal{H}_{\nu}^{2}(T_{\Omega})$.
\end{theorem}

\begin{proof}
 Denote temporarily $\mathcal L^\Ome=\mathcal L_\nu^\Ome$.
We define $\mathcal L^\Ome$
on the domain
$$
\{f\in L^2_\nu(\Ome): \int_\Ome e^{-(\cdot|x)}
f(x)d\mu_{\nu}(x)\in L^2_\nu(\Ome)
\}.
$$
Note that by Cauchy-Schwarz inequality the condition $f\in
L^2_\nu(\Ome)$ implies that $\mathcal L^\Ome f$ is a well-defined
function on $\Ome$ since $e^{-(y|x)}$ for fixed $y\in \Ome$ is a
function in $L^2_\nu(\Ome)$. It is then easy to prove that
$\mathcal L^\Ome$ so defined is a self-adjoint positive operator.
Let $f\in L^2_\nu(\Ome)$ be in the domain of $RR^{\ast}$, then
\begin{equation}
\label{RR-ast}
\begin{split}
RR^{\ast}f(y)
&  =R^{\ast}f(y)\\
&  =(R^{\ast}f\mid k_{y})_{\mathcal{H}_{\nu}(T(\Omega))}\\
&  =(f\mid Rk_{y})_{L_{\nu}^{2}(\Omega)}\\
&  =\int_{\Omega}f(x)\overline{K(x,y)}\,d\mu_{\nu}(x)\\
&  =\Gamma_{\Omega}(\nu)\int_{\Omega}f(x)\Delta(x+y)^{-\nu}\,d\mu_{\nu}(x)\\
&  =\int_{\Omega}f(x)\mathcal{L}^\Ome(e^{-(y\mid\cdot)})(x)\,d\mu_{\nu}(x)\\
&  =\int_{\Omega}e^{-(y\mid x)}\Delta(x)^{\nu-d/r}\mathcal{L}^\Ome(f)(x)\,dx\\
&={\mathcal L^\Ome}
(\mathcal L^\Ome f)(y),
\end{split}
\end{equation}
and
$$(\mathcal L^\Ome f,
\mathcal L^\Ome f)=(R^{\ast}f, R^{\ast}f)<\infty.
$$
Thus  $f$ is in the domain of $(\mathcal L^\Ome)^2$ and
$RR^{\ast}=(\mathcal L^\Ome)^2$, since $(\mathcal L^\Ome)^2$ is a
self-adjoint extension of $RR^\ast$, which is also self-adjoint by
the von Neumann theorem (see e.g. \cite{RS-MMMPH},  VIII. problem
45.)

Consider the inverse operator $R^{-1}$
acting on the image of $R$. For a function $g$ in
the image of $R$,  $R^{-1} g$ is the unique
extension of $g$ to a holomorphic function on $T(\Ome)$.
Thus
$R^{-1}{\mathcal L^\Ome}={\mathcal L}_\nu$.
Let $R^{-1}$ act on the previous equality
(\ref{RR-ast})
$$
R^\ast f={\mathcal L}_\nu {\mathcal L^\Ome} f.
$$
This proves the polar decomposition formula. Since $R^\ast
$ is densily defined and $R$ is an injective closed
operator we have that the unitary
part ${\mathcal L}_\nu$
extends to a unitary operator.
\end{proof}
\begin{remark} The multiplication map $f(x)\rightarrow \Delta(x)^{\frac{\nu}{2}}f(x)$
induces a unitary isomorphism between $L_\nu ^2(\Omega)$ to
$L^2(\Omega,d\mu_\circ)$. In turn, $L^2(\Omega,d\mu_\circ)$ is
unitarily equivalent to $L^{2}(\mathcal{D}_{\mathbb{R}},d\eta)$
via a scalar multiple of the Cayley transform.  This is likewise
true for $\mathcal{H}_{\nu}(T(\Omega))$ and $\mathcal{H}_{\nu
}(\mathcal{D})$ (see Lemma \ref{le-HnuD}).  Furthermore, these
isomorphisms intertwine the corresponding restriction maps.  It
follows then from Theorem \ref{th-main} that $R$ is a continuous
operator for $\nu>a(r-1)$.

\end{remark}

\begin{corollary}
Assume that $\nu>1+a(r-1)$. Then $\mathcal{L}_{\nu}^{\ast}:\mathcal{H}_{\nu
}^{2}(T(\Omega))\rightarrow L_{\nu}^{2}(\Omega)$ is given by the integral
operator
\[
\mathcal{L}_{\nu}^{\ast}F(x)=\alpha_{\nu}\int_{T(\Omega)}F(z)e^{-(\overline
{z}\mid x)}\Delta(y)^{\nu-2d/r}\,dx\,dy
\]
where $\alpha_{\nu}=\frac{2^{r\nu}}{(4\pi)^{d}\Gamma_{\Omega}(\nu-d/r)}$.
\end{corollary}

\begin{proof}
Write $\mathcal{L}$ for $\mathcal{L}_{\nu}$. Let $F\in\mathcal{H}_{\nu}%
^{2}(T(\Omega))$ and $f\in L_{\nu}^{2}(\Omega)$. Then
\begin{align*}
(\mathcal{L}^{\ast}F   \mid
f)&=(F\mid\mathcal{L}f)_{\mathcal{H}_{\nu
}^2(T(\Omega))}\\
&  =\alpha_{\nu}\int_{T(\Omega)}F(z)\overline{\mathcal{L}f(z)}\Delta
(y)^{\nu-2d/r}\,dx\,dy\\
&  =\alpha_{\nu}\int_{T(\Omega)}F(z)\int_{\Omega}\overline{e^{-(z\mid t)}%
f(t)}\Delta(t)^{\nu-d/r}\Delta(y)^{\nu-2d/r}\,dt\,dx\,dy\\
&  =\alpha_{\nu}\int_{\Omega}\int_{T(\Omega)}F(z)e^{-(\overline{z}\mid
t)}\Delta(y)^{\nu-2d/r}\,dx\,dy\,\overline{f(t)}\,d\mu_{\nu}(t).
\end{align*}
\end{proof}

Using the Laplace transform we transfer the $G^{\gamma}$-action,
$\pi_{\nu}$,
on $\mathcal{H}_{\nu}^{2}(T(\Omega))$ to an equivalent action, denoted by
$\lambda_{\nu}$, on $L_{\nu}^{2}(\Omega)$ and notice the following simple fact:

\begin{lemma}
\label{le-Haction}The $H^{\gamma}$-action on $L_{\nu}^{2}(\Omega)$ is given by
the formula
\[
\lambda_{\nu}(h)f(x)=\mathrm{Det}(h)^{\frac{\nu}{p}}f(h^{t}x),
\]
where the determinant is taken as a real linear transformation on $J$ and
$h^{t}$ denotes the transpose with respect to the real form $(x\mid y).$
\end{lemma}

\begin{proof}
For $h\in H^{\gamma}$ a straightforward calculation gives that
$J(h,z)=\mathrm{Det}(h)$ and $d\mu_{\nu}(hx)=\mathrm{Det}(h)^{\frac{\nu r}{d}%
}d\mu_{\nu}(x)$. Thus, for $f\in L_{\nu}^{2}(\Omega)$ we have
\begin{align*}
\pi_{\nu}(h)\mathcal{L}_{\nu}f(z)  &  =J(h^{-1},z)^{\frac{\nu r}{2d}}%
\int_{\Omega}e^{-(h^{-1}z|x)}f(x)d\mu_{\nu}(x)\\
&  =\mathrm{Det}(h)^{-\frac{\nu r}{2d}}\int_{\Omega}e^{-(z\mid x)}%
f(h^{t}x)d\mu_{\nu}(h^{t}x)\\
&  =\mathrm{Det}(h)^{\frac{\nu r}{2d}}\mathcal{L}_{\nu}(f\circ h^{t})(z).
\end{align*}
This calculation now implies the lemma.
\end{proof}

We follow \cite{FK-book}, p. 343, and define the generalized Laguerre
polynomials by the formula
\[
L_{\mathbf{m}}^{\nu}(x)=(\nu)_{\mathbf{m}}\sum_{|\mathbf{n|\leq|m|}}%
\begin{pmatrix}
\mathbf{m}\\
\mathbf{n}%
\end{pmatrix}
\frac{1}{(\nu)_{\mathbf{n}}}\psi_{\mathbf{n}}(-x)\,,
\]
and the generalized Laguerre functions by
\[
\ell_{\mathbf{m}}^{\nu}(x)=e^{-\mathrm{Tr}(x)}L_{\mathbf{m}}^{\nu}(2x)\,.
\]
By Proposition XV.4.2 in \cite{FK-book}, p. 344, we get:.

\begin{theorem}
The Laguerre functions form an orthogonal basis of $L_{\nu}^{2}(\Omega)^{L}$.
Furthermore,
\[
\mathcal{L}_{\nu}(\ell_{\mathbf{m}}^{\nu})=\Gamma_{\Omega}(\mathbf{m}%
+\nu)q_{\mathbf{m}}^{\nu}\,.
\]
\end{theorem}

Let $E$ be the Euler operator on $\Omega$ (or $V$). Specifically,
\[
Ef(x)=\frac{d}{dt}f(tx)|_{t=1}=\frac{d\,}{dt}f(\operatorname{exp}(tZ_{0})\cdot
x)|_{t=0}\,.
\]
We now obtain the following recursion formula.

\begin{theorem}
\label{The-Euler}The Laguerre functions are related by the following recursion
relations
\[
2E\ell_{\mathbf{m}}^{\nu}=-\nu r\ell_{\mathbf{m}}^{\nu}-\sum_{j=1}^{r}%
\begin{pmatrix}
\mathbf{m}\\
\mathbf{m-\gamma_{j}}%
\end{pmatrix}
(m_{j}-1+\nu-(j-1)\frac{d}{2})\ell_{\mathbf{m-\gamma_{j}}}^{\nu}+\sum
_{j=1}^{r}c_{\mathbf{m}}(j)\ell_{\mathbf{m+\gamma_{j}}}^{\nu}.
\]
\end{theorem}

\begin{proof}
By Lemma \ref{le-Haction} the $H^{\gamma}$-action is given by
\[
\lambda_{\nu}(h)f(x)=\operatorname{det}(h)^{\frac{\nu r}{2d}}f(h^{t}x),
\]
The infinitesimal action is then given in the usual way by
differentiation. Thus by Lemma \ref{le-ggamma}
$Z_{0}\in\mathfrak{h}^{\gamma}$ and $Z_{0}$ acts on the smooth
vectors in $L_{\nu}^{2}(\Omega)$ by :
\[
\lambda_{\nu}(Z_{0})f(x)=(\frac{\nu r}{2}+E)f(x).
\]
According to Lemma \ref{le-rerelqmnu}, part 2, we have%

\[
-\pi_{\nu}(2Z_{0})q_{\mathbf{m},\nu}=\sum_{j=1}^{r}%
\begin{pmatrix}
\mathbf{m}\\
\mathbf{m-\gamma_{j}}%
\end{pmatrix}
q_{\mathbf{m-\gamma_{j},\nu}}-\sum_{j=1}^{r}(\nu+\mathbf{m}_{j}-\frac{a}%
{2}(j-1)c_{\mathbf{m}}(j))q_{\mathbf{m+\gamma_{j},\nu}}.
\]
The inverse Laplace transform of the above equation gives
\begin{align*}
-2(\frac{r\nu}{2}+E)\frac{\ell_{\mathbf{m}}^{\nu}}{\Gamma(\nu+\mathbf{m})}
&
=\sum_{j=1}^{r}%
\begin{pmatrix}
\mathbf{m}\\
\mathbf{m-\gamma_{j}}%
\end{pmatrix}
\frac{\ell_{\mathbf{m-\gamma_{j}}}^{\nu}}{\Gamma(\nu+\mathbf{m-\gamma_{j}})}\\
\,  &  -\sum_{j=1}^{r}(\nu+\mathbf{m}_{j}-\frac{d}{2}(j-1)c_{\mathbf{m}%
}(j))\frac{l_{\mathbf{m+\gamma_{j}}}^{\nu}}{\Gamma(\nu+\mathbf{m+\gamma_{j}}%
)}.
\end{align*}
This simplifies to
\[
-2(\frac{\nu r}{2}+E)l_{\mathbf{m}}^{\nu}=\sum_{j=1}^{r}%
\begin{pmatrix}
\mathbf{m}\\
\mathbf{m-\gamma_{j}}%
\end{pmatrix}
(m_{j}-1+\nu-(j-1)\frac{d}{2})\ell_{\mathbf{m-\gamma_{j}}}^{\nu}-\sum
_{j=1}^{r}c_{\mathbf{m}}(j)\ell_{\mathbf{m+\gamma_{j}}}^{\nu},
\]
and proves the theorem.
\end{proof}

In the above Theorem we used the action of $Z_{0}$, but to
directly derive a differential equation satisfied by the Laguerre
functions one uses the element $\xi$ and part 1 in Lemma
\ref{le-rerelqmnu}:

\begin{remark}
If the Laguerre polynomials were defined by the formula
\[
^{\circ}L_{\mathbf{m}}^{\nu}(x)=\frac{(\nu)_{\mathbf{m}}}{\Gamma
(\nu+\mathbf{m})}\sum_{|\mathbf{n|\leq|m|}}%
\begin{pmatrix}
\mathbf{m}\\
\mathbf{n}%
\end{pmatrix}
\frac{1}{(\nu)_{\mathbf{n}}}\psi_{\mathbf{n}}(-x)=\frac{L_{\mathbf{m}}^{\nu
}(x)}{\Gamma(\nu+\mathbf{m})},
\]
it would agree with our definition given in \cite{doz1}. With this definition
the formula in the proceeding theorem would become
\[
2E{^{\circ}\ell_{\mathbf{m}}^{\nu}}=-\nu r^{\circ}\ell_{\mathbf{m}}^{\nu}%
-\sum_{j=1}^{r}%
\begin{pmatrix}
\mathbf{m}\\
\mathbf{m-\gamma_{j}}%
\end{pmatrix}
{^{\circ}\ell_{\mathbf{m-\gamma_{j}}}^{\nu}}+\sum_{j=1}^{r}(m_{j}%
-1+\nu-(j-1)\frac{d}{2})c_{\mathbf{m}}(j){^{\circ}\ell_{\mathbf{m+\gamma_{j}}%
}^{\nu}}.
\]
\end{remark}

\begin{remark}
Theorem \ref{The-Euler} involves both a creation and an
annihilation operator, i.e., it involves both a step up
$\mathbf{n}\mapsto\mathbf{n}+\gamma_{j}$ and a step down
$\mathbf{n}\mapsto\mathbf{n}-\gamma_{j}$. This is related to the
fact that the element $Z_{0}\in\mathfrak{z}(\mathfrak{h})$ which
is used to derive the relation in Theorem \ref{The-Euler} has a
decomposition $\xi=E_{+}+E_{-}$ into $H^{\gamma}$ invariant
element where $E_{+}$ is in the $\mathfrak{p}^{\gamma+}$ and steps
down and $E_{-}$ is in $\mathfrak{p}^{\gamma-}$ and steps up. The
elements $E_{+}$ and $E_{-}$ are not in
$\mathfrak{h}_{\mathbb{C}}^{\gamma}$ and hence they act as a
\textbf{second order} differential operator. In the case of the
upper half-plane (in particular, see Proposition 2.7 and Theorem
3.4 in \cite{doz1} where we are using $\mathbb{R}+i\mathbb{R}^{+}$
as
a realization of the tube domain) this corresponds to%
\[
\left(
\begin{array}
[c]{cc}%
0 & 1\\
1 & 0
\end{array}
\right)  =\frac{1}{2}\left(
\begin{array}
[c]{cc}%
-i & 1\\
1 & i
\end{array}
\right)  +\frac{1}{2}\left(
\begin{array}
[c]{cc}%
i & 1\\
1 & -1
\end{array}
\right)
\]
and%
\begin{align*}
\frac{1}{2}\left(
\begin{array}
[c]{cc}%
-i & 1\\
1 & i
\end{array}
\right)   &  \leftrightarrow\frac{-i}{2}\left(  tD^{2}+(2t+\nu\right)
D+(t+\nu)\\
\frac{1}{2}\left(
\begin{array}
[c]{cc}%
i & 1\\
1 & -1
\end{array}
\right)   &  \leftrightarrow\frac{-i}{2}\left(  tD^{2}-(2t-\nu\right)
D+(t-\nu)\,.
\end{align*}
on $i\mathbb{R}^{+}$. Furthermore, one can in a similar way get a
differential equation satisfied by the Laguerre functions by
applying the Cayley transform
to part 1 in Lemma \ref{le-rerelqmnu}%
\[
\pi_{\nu}(\xi)p_{\mathbf{m},\nu}=(r\nu+2\left|  \mathbf{m}\right|
)p_{\mathbf{m},\nu}\,.
\]
This equation is also a second order differential equation, which
in the case
of the upper half-plane is%
\[
(tD^{2}+\nu D-t)\ell_{n}^{\nu}=-(\nu+2n)\ell_{n}^{\nu}\,.
\]
All of these equations can be caried over to the general case to
find the radial part of the corresponding operators. We think that
it is an interesting problem to find an explicit formula for the
raising operator, annihilator operator and the operator
$\pi_{\nu}(\xi)$, in the general case.
\end{remark}

\providecommand{\bysame}{\leavevmode\hbox to3em{\hrulefill}\thinspace}

%
%
%
%
%
%
%
%
%
%
%
%
%
%
%
%
%
\end{document}